\documentclass[11pt]{article}
\usepackage{amsfonts, amsmath, amssymb, enumitem}

\PassOptionsToPackage{obeyspaces}{url}
\usepackage[colorlinks=true,citecolor=blue,urlcolor=blue,linkcolor=blue,bookmarksopen=true]{hyperref}
\usepackage{amsthm}
\usepackage{tikz}
\usetikzlibrary{calc,shapes.geometric,positioning,decorations.pathmorphing}
\usepackage{mathtools}

\usepackage{longtable}

\usepackage{etoolbox}
\patchcmd{\thebibliography}{\leftmargin\labelwidth}{\leftmargin\labelwidth\addtolength\itemsep{-0.1\baselineskip}}{}{}

\author{Boris Bukh\thanks{Department of Mathematical Sciences, Carnegie Mellon University, Pittsburgh, PA, USA. \texttt{bbukh@math.cmu.edu}. Supported in part by Sloan\
	Research Fellowship and by U.S.\ taxpayers through NSF CAREER grant DMS-1555149.}\and
	Christopher Cox\thanks{Department of Mathematical Sciences, Carnegie Mellon University, Pittsburgh, PA, USA. \texttt{cocox@andrew.cmu.edu}. Supported in part by U.S.\ taxpayers through NSF CAREER grant DMS-1555149.}}

\title{Periodic words, common subsequences and frogs}
\date{}

\usepackage[nameinlink]{cleveref}

\newtheorem{theorem}{Theorem}

\newtheorem{lemma}[theorem]{Lemma}

\newtheorem{proposition}[theorem]{Proposition}
\crefname{proposition}{proposition}{propositions}

\newtheorem{defn}[theorem]{Definition}
\crefname{defn}{definition}{definitions}

\newtheorem{claim}[theorem]{Claim}
\crefname{claim}{claim}{claims}

\newtheorem{conj}[theorem]{Conjecture}
\crefname{conj}{conjecture}{conjectures}

\crefname{remark}{remark}{remarks}

\crefname{enumi}{part}{parts}

\let\theparentequation\theequation
\patchcmd{\theparentequation}{equation}{parentequation}{}{}

\newcommand*{\eqdef}{\stackrel{\mbox{\normalfont\tiny def}}{=}}   % definition by equality
\newcommand*{\abs}[1]{\lvert #1\rvert}                           % Absolute values, cardinality
\newcommand*{\N}{\mathbb{N}}                                     % Natural numbers
\newcommand*{\Z}{\mathbb{Z}}                                     % Integers
\newcommand*{\R}{\mathbb{R}}

\newcommand*{\Rand}[1]{R_{#1}}
\renewcommand*{\hat}[1]{\widehat{#1}}
\renewcommand*{\cal}[1]{\mathcal{#1}}
\newcommand*{\tod}{\stackrel{d}{\to}}
\DeclareMathOperator{\E}{\mathbb{E}}                             % Expectation
\renewcommand*{\Pr}{\mathbf{Pr}}                           % Probability
\DeclareMathOperator{\Var}{Var}                                  % Variance
\DeclareMathOperator{\len}{len}                                  % Length of a word
\DeclareMathOperator{\LCS}{LCS}                                  % Length of the longest common subsequence
\def\mydot at (#1,#2){\fill (#1,#2) circle [radius=0.2]}         % Dots in profile pictures
\newcommand*{\beg}{\mathtt{begin}}
\newcommand*{\nd}{\mathtt{end}}
\newcommand*{\trans}{\mathtt{trans}}
\newcommand*{\Mod}[1]{\ (\mathrm{mod}\ #1)}						% nicer mod
\newcommand*{\bet}{\Sigma}										% alphabet
\newcommand*{\F}{\widetilde{F}}                                 % for the intermediate frog arrangements

\def\froggie{\protect\includegraphics[scale=0.014]{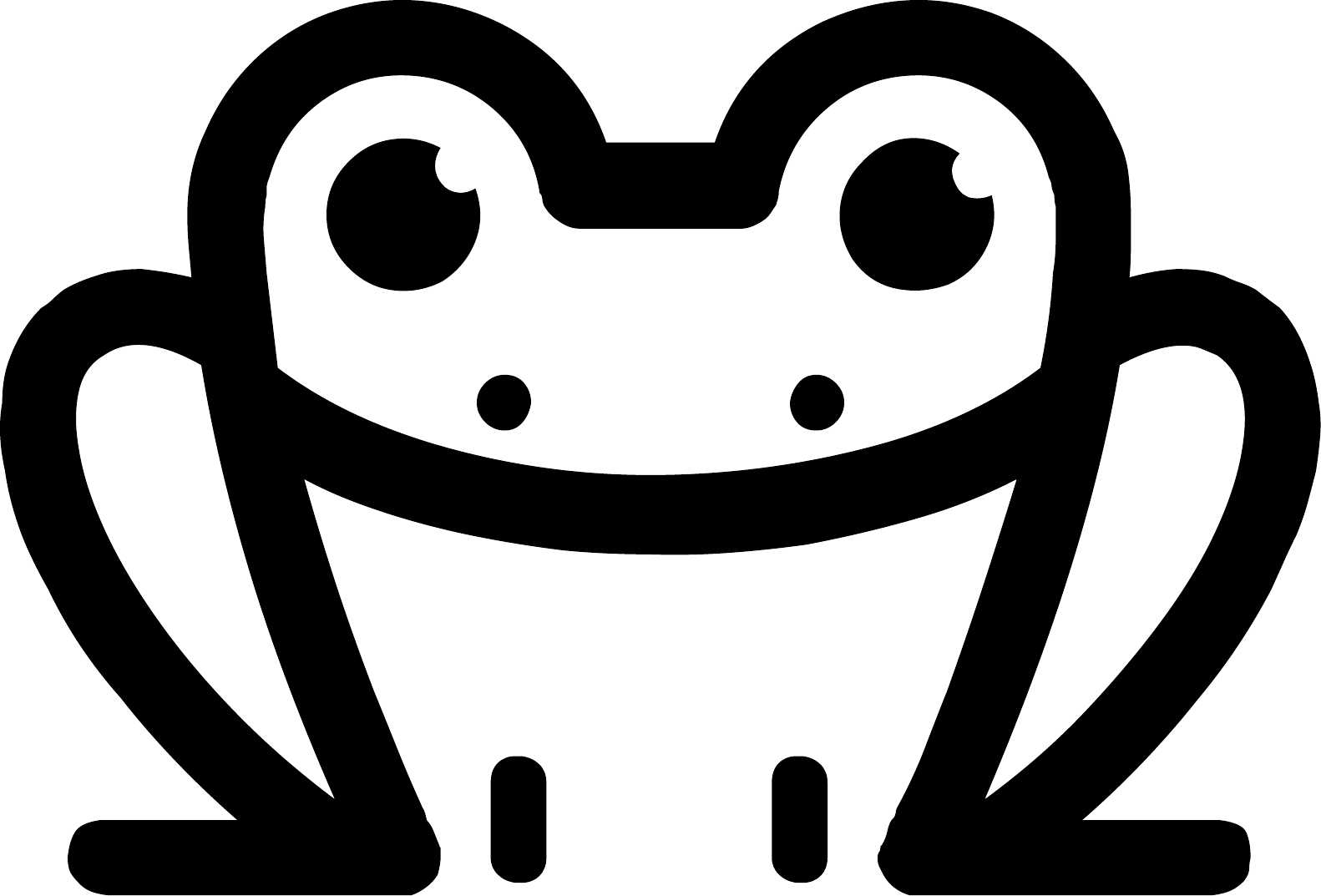}}
\def\lilypad{\raisebox{-0.3ex}{\includegraphics[height=0.66em]{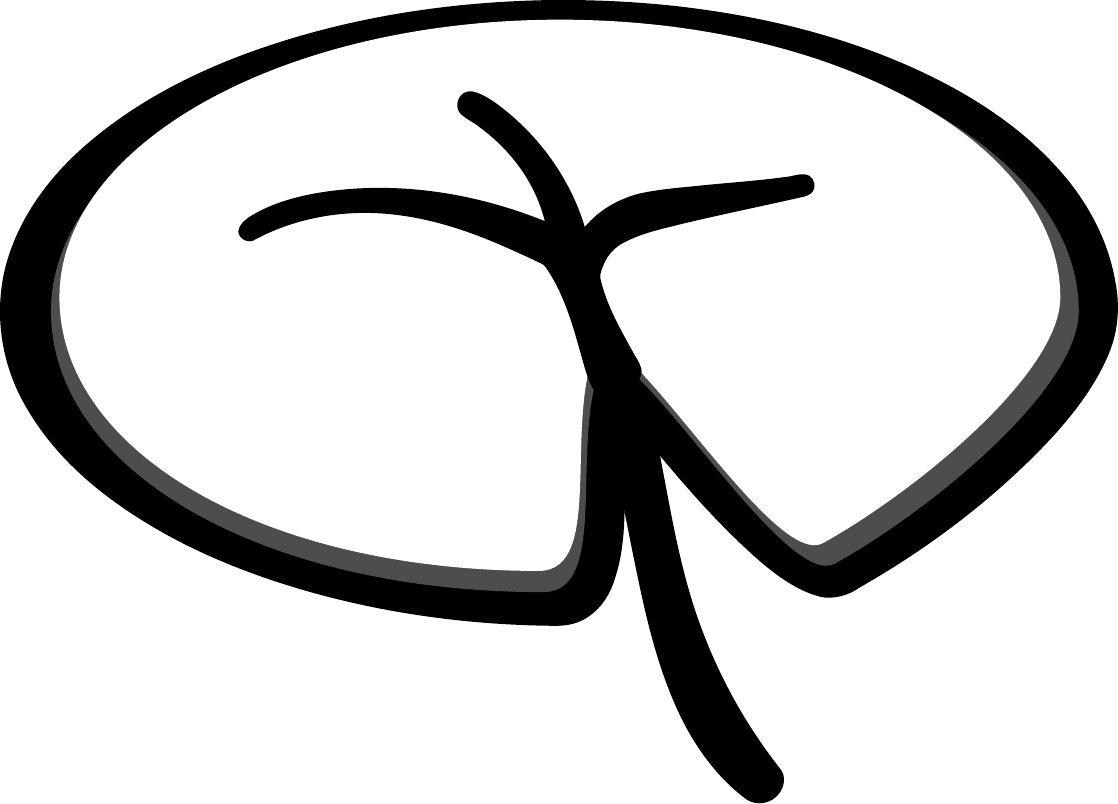}}}

\begin{document}

\maketitle

\begin{abstract}
	Let $W^{(n)}$ be the $n$-letter word obtained by repeating a fixed word~$W$, and let $\Rand{n}$
	be a random $n$-letter word over the same alphabet. We show several results about
	the length of the longest common subsequence (LCS) between $W^{(n)}$ and $\Rand{n}$; in particular,
	we show that its expectation is $\gamma_W n-O(\sqrt{n})$ for an efficiently-computable constant~$\gamma_W$.

	This is done by relating the problem to a new interacting particle system, which we dub ``frog dynamics''.
	In this system, the particles (`frogs') hop over one another in the order given by their labels.
	Stripped of the labeling, the frog dynamics reduces to a variant of the PushTASEP.

	In the special case when all symbols of $W$ are distinct, we obtain an explicit formula for the constant~$\gamma_W$ and a closed-form expression for the stationary distribution of the associated frog dynamics.

	In addition, we propose new conjectures about the asymptotic of the LCS of a pair of random words.
	These conjectures are informed by computer experiments using a new heuristic algorithm to compute the LCS.
        Through our computations, we found periodic words that are more random-like than a random word, as measured by the LCS.
\end{abstract}

\section{Introduction}
\paragraph{The longest common subsequence problem.}
A \emph{word} is a finite sequence of symbols from some alphabet. We denote by $\len W$
the length of the word~$W$. A \emph{subsequence} of a word $W$ is a word obtained from
$W$ by deleting some symbols from~$W$; the symbols in a subsequence are not required to appear contiguously in $W$. A \emph{common subsequence} between words $W$ and $W'$
is a subsequence of both $W$ and~$W'$. We denote by $\LCS(W,W')$ the length of the
\emph{longest common subsequence} between $W$ and~$W'$. We write $W_i$ for the $i$'th symbol of $W$,
with indexing starting from~$0$.

Throughout the paper, we use $\bet$ to denote the alphabet, and we
write $R\sim \bet^n$ to indicate that $R$ is a word chosen uniformly at random
from~$\bet^n$. A long-standing problem is to understand $\LCS(R,R')$
for a pair of independently chosen words $R,R'\sim \bet^n$.
Whereas it is known that
\begin{equation}\label{eq:chvatal_sankoff}
	\E \LCS(R,R')=\gamma n+o(n)
\end{equation}
for some constant $\gamma$ depending on $\abs{\bet}$, little else is known.
We mention three open problems.
\begin{enumerate}
	\item The rate of convergence in \eqref{eq:chvatal_sankoff} is unknown. The original
		proof of \eqref{eq:chvatal_sankoff} by Chv\'atal and Sankoff~\cite{chvatal_1975} did not supply
		any bound on the $o(n)$ term.
                Alexander~\cite{alexander_1994} showed that
		$\E \LCS(R,R')=\gamma n+O(\sqrt{n \log n})$.
	\item The value of $\gamma$, which is often called the \emph{Chv\'atal--Sankoff constant},
		is unknown. The best rigorous bounds for the binary alphabet are due to Lueker~\cite{lueker_2009},
		whereas Kiwi, Loebl and Matou\v{s}ek~\cite{kiwi_2005} gave an asymptotic for $\gamma$ as $|\bet|\to\infty$.
	\item It is believed that $\LCS(R,R')$ is approximately normal, and that its variance is linear in~$n$.
		Yet it is not even known that $\Var \LCS(R,R')$ tends to infinity with~$n$.
\end{enumerate}

We performed extensive computer simulations using a new heuristic algorithm in order to compute $\LCS(R,R')$ for large $n$.
These simulations suggest that $\E\LCS(R,R')=\gamma n-\Theta(n^{1/3})$ and that $\gamma\approx 0.8122$ for the binary alphabet.
We shall discuss both the algorithm and the computer simulations in \Cref{sec:computerjoin}.

\paragraph{Periodic words.}
A word $W$ is \emph{$k$-periodic} if $W_{i+k}=W_i$ holds for all values of $i$, for which
both sides are defined (that is for $i=0,1,\dotsc,\len W-k-1$).
For a word $W$ of length $k$, write $W^{(n)}$ for the $k$-periodic word
of length $n$ which is obtained by repeating $W$ the appropriate number of times
(which might be fractional if $k$ does not divide~$n$). For example,
if $W=aba$, then $W^{(8)}=abaabaab$.
Additionally, write $W^{(\infty)}$ to denote the $k$-periodic word obtained by repeating $W$ ad infinitum.

\begin{figure}[ht]
	\begin{center}
		\begin{tikzpicture}[scale=4]
	\node[above left] (gamma) at (0,1.1) {$\gamma_W(\rho)$};
	\node[below right] (rho) at (3,0) {$\rho$};
	\draw (-0.1,0)--(3,0);
	\draw (0,-0.1)--(0,1.1);
	\node[below] (s1x) at (1/4,0) {$s_1$};
	\node[below] (s2x) at (5/12,0) {$s_2$};
	\node[below] (s3x) at (5/6,0) {$s_3$};
	\node[below] (s4x) at (5/2,0) {$s_4$};
	\node[left] (s1y) at (0,1/4) {$1/4$};
	\node[left] (s2y) at (0,3/8) {$3/8$};
	\node[left] (s3y) at (0,7/12) {$7/12$};
	\node[left] (s4y) at (0,1) {$1$};
	\node[circle,fill=black] (s1xy) at (1/4,1/4) {};
	\node[circle,fill=black] (s2xy) at (5/12,3/8) {};
	\node[circle,fill=black] (s3xy) at (5/6,7/12) {};
	\node[circle,fill=black] (s4xy) at (5/2,1) {};
	\draw[dotted] (s1x)--(s1xy)--(s1y);
	\draw[dotted] (s2x)--(s2xy)--(s2y);
	\draw[dotted] (s3x)--(s3xy)--(s3y);
	\draw[dotted] (s4x)--(s4xy)--(s4y);
	\draw[ultra thick] (0,0)--(s1xy)--(s2xy)--(s3xy)--(s4xy)--(3,1);
\end{tikzpicture}
	\end{center}
	\caption{\label{fig:gamma} The plot of $\gamma_W(\rho)$ for $\bet=[4]$ and $W=1234$. Here, $s_1=1/4$, $s_2=5/12$, $s_3=5/6$ and $s_4=5/2$. Furthermore, $\tau_W(s_1)=\sqrt{{3\over 512\pi}}$, $\tau_W(s_2)=\sqrt{{145\over 13824\pi}}$, $\tau_W(s_3)=\sqrt{{79\over 3456\pi}}$, $\tau_W(s_4)=\sqrt{{5\over 128\pi}}$ and $\tau_W(\rho)=0$ otherwise.}
\end{figure}
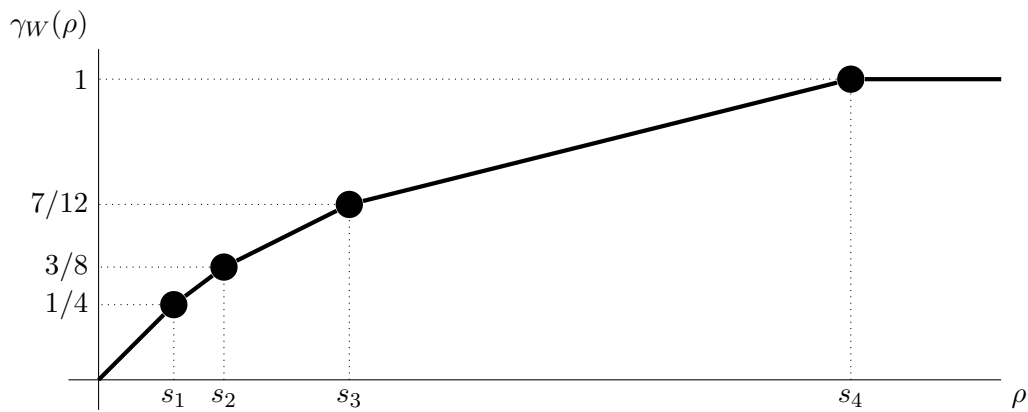

In attempt to better understand the problems enumerated above, in this paper we tackle a simpler
random variable $\LCS(R,W^{(n)})$ where $W$ is a fixed word. This random variable was previously studied by Matzinger--Lember--Durringer \cite{matzinger_variance}.
We give answers to the analogues of all three problems that we stated above.
These answers are summarized in the following theorem; more precise results are below in \Cref{thm:frog}.
For a visualization of the following theorem, see \Cref{fig:gamma}.

\begin{theorem}\label{thm:basic}
	Let $\rho$ be a positive real number. Fix $W\in \bet^k$ and let $R\sim \bet^n$ be an $n$-letter random word. Then
	\[
		\E \LCS(R,W^{(\rho n)})=\gamma_W n-\tau_W\sqrt{n}+O(1),
	\]
	where
	\begin{enumerate}[label=(\roman*)]
		\item \label{part:piecewiselinear} $\gamma_W=\gamma_W(\rho)$ is a non-negative piecewise linear function of~$\rho$.
        \item \label{part:slope} The slope of $\gamma_W(\rho)$ is a non-increasing function of $\rho$.
		\item \label{part:positive} $\tau_W=\tau_W(\rho)$ is nonzero only at the points where the slope of $\gamma_W(\rho)$ changes,
			and $\tau_W$ is strictly positive at those points.
		\item \label{part:normality} The random variable $\LCS(R,W^{(\rho n)})$ is asymptotically normal with linear variance if $\rho>1/|\bet|$, $\tau_W(\rho)=0$ and either
			\begin{enumerate}
				\item the slope of $\gamma_W(\rho)$ is positive, or
				\item there is some symbol in $\bet$ which does not appear in $W$.
			\end{enumerate}

            If $\tau_W(\rho)\neq 0$, then $\LCS(R,W^{(\rho n)})$ still has linear variance but is \emph{not} asymptotically normal.

            In all other cases, $\LCS(R,W^{(\rho n)})$ has sub-linear variance.
		\item\label{part:algorithm} There exists an algorithm that computes $\gamma_W$ and $\tau_W$ from~$W$.
	\end{enumerate}
\end{theorem}

From \cref{part:positive}, it is clear that $\tau_W\neq 0$ happens rarely. However, it does happen
for infinitely many $W$ even in the case $\rho=1$; see \Cref{thm:nice} for examples.

\Cref{part:normality} extends a result of Matzinger--Lember--Durringer~\cite{matzinger_variance}, who showed that $\Var\LCS(R,W^{(n)})$ is linear in $n$ when $|\bet|=2$.
The most interesting piece of \cref{part:normality} is that $\LCS(R,W^{(\rho n)})$ is \emph{not always} normal.
In fact, when $\tau_W(\rho)\neq 0$, we show that, under the correct shifting and scaling, this variable converges to the minimum of two Gaussians (see \cref{nnn:hitspeed} of \Cref{thm:normnotnorm} for the precise statement).

\Cref{part:algorithm} contrasts the constants $\gamma_W$ with the usual Chv\'atal--Sankoff constant $\gamma$ from \eqref{eq:chvatal_sankoff}. Whereas convergence of limits defining both
$\gamma$ and $\gamma_W$ is an easy application of a standard superadditivity argument, $\gamma$ is not known to be computable in finite time. Even the algorithms to approximate
$\gamma$ are non-trivial, see \cite{lueker_2009}.

\begin{figure}[ht]
	\centering
	\input{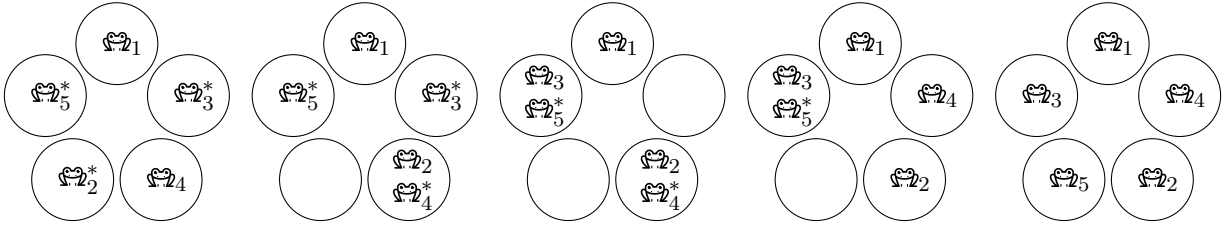}
	\caption{The sequence of frog hops resulting from poking frogs $\froggie_2$, $\froggie_3$ and $\froggie_5$. Here, frogs move in the anti-clockwise direction, and a ${}^*$ indicates that the frog is agitated.}
	\label{fig:frogex}
\end{figure}

\paragraph{Frog dynamics.} The key to \Cref{thm:basic} is the analysis of the following dynamical system.
Let $W$ be a fixed word, and set $k=\len W$. Imagine a circle of $k$ lily pads, each of which is occupied by a frog. The $k$ frogs vary
from a large nasty frog to a little harmless froggie. No two frogs are equally nasty, and are thus linearly ordered by their nastiness. They all face in the same (circular) direction.
At each time step $t=0,1,\dotsc$, the following happens:
\begin{enumerate}
	\item The monster living below pokes some of the frogs with its tentacles. Each poked frog gets agitated, and wants to jump away.
	\item In the order of descending nastiness, starting from the nastiest frog, each of the agitated frogs will leap
		to the next `available' lily pad, that is either empty or occupied by a less menacing frog.
		Doing so causes the current occupant to become agitated, and the frog that just hopped calms down.

		This process repeats until all frogs are content once more.
\end{enumerate}
Note that, with each step of this process, the nastiest agitated frog get less and less nasty. This
guarantees termination of the process, and implies that no frog jumps over another agitated frog.
Below, in \Cref{fig:frogex}, is an example of one round of this process, where here and thereafter we denote
the frogs $\froggie_1, \froggie_2,\dotsc,\froggie_k$ in the order of nastiness, with $\froggie_1$ being the nastiest.

We denote the lily pads $\lilypad_0,\dotsc,\lilypad_{k-1}$ in the circular order.
The lily pads correspond to symbols of $W$, and so we label $\lilypad_i$ with~$W_i$ (recall that we index symbols from $0$).

We record the frogs' positions in a \emph{frog arrangement}, which formally is just a bijection from $\{\froggie_1,\dotsc,\froggie_k\}$ to~$\{\lilypad_0,\dotsc,\lilypad_{k-1}\}$. We denote by $\mathcal{F}$ the collection
of all frog arrangements. For a frog arrangement $F\in \mathcal{F}$ and a symbol $a\in \bet$, we let $Fa$ be the frog
arrangement resulting from poking all lily pads labeled~$a$, and waiting for the ensuing frenzy to settle.
For a word $R=R_0\dotsb R_{\ell-1}$, we write $FR\eqdef FR_0\dotsb R_{\ell-1}$. Note that this notation respects
concatenation, i.e., $(FR)R'=F(RR')$ for any two words~$R,R'$.
Denote by $D_m(F,R)$ the total displacement of $\froggie_m$ as the word $R$ is applied to the frog arrangement $F$.
\medskip

\noindent\textbf{Example.} The frog arrangement shown in the left-most image in \Cref{fig:frogex} has $F(\froggie_1)=\lilypad_0$, $F(\froggie_2)=\lilypad_2$, $F(\froggie_3)=\lilypad_4$, $F(\froggie_4)=\lilypad_3$ and $F(\froggie_5)=\lilypad_1$.
If $W=abbab$, then the right-most image in \Cref{fig:frogex} is precisely the arrangement $Fb$ and $D_1(F,b)=0$, $D_2(F,b)=D_4(F,b)=D_5(F,b)=1$ and $D_3(F,b)=2$.
\medskip

\noindent\textbf{Example.} Suppose that $F$ is a frog arrangement wherein the label of $F(\froggie_k)$ is $a$ and the symbol $a$ appears only once in $W$.
By poking lily pads labeled $a$, only $\froggie_k$ becomes agitated and it will hop all the way around the circle of lily pads back to its original pad.
Therefore, $Fa=F$ with $D_1(F,a)=\dots=D_{k-1}(F,a)=0$ and $D_k(F,a)=k$.
\medskip

Starting with a frog arrangement $F_0$, set $F_i=F_0R_0\dotsb R_{i-1}$, where the symbols $R_0,R_1,\dotsc$ are chosen independently at random from $\bet$. Since $F_{i+1}=F_iR_i$,
the sequence $F_0,F_1,\dotsc$ forms a Markov chain. We call this Markov chain \emph{the frog dynamics} associated with~$W$. (We note that
our frog dynamics has no mathematical relation to the model of simple random walks on $\Z^d$ known as the `frog model' and studied for example
in \cite{notfrogs}.)

Observe that if $W$ can be written as, say, $W=UU$ for some other word $U$, then $W^{(n)}=U^{(n)}$. Hence in this case we may as well use the shorter word $U$ in lieu of~$W$.
To this end, we say that a word $W\in \bet^k$ is \emph{reducible} if there is some other word $U\in \bet^\ell$ with $W=U^{(k)}$ where $\ell<k$ and $\ell\mid k$.
Otherwise, we say that $W$ is \emph{irreducible}.

\begin{theorem}\label{thm:frog}
	Let $F_0,F_1,\dotsc$ be the frog dynamics associated with an irreducible word $W\in \bet^k$. Then
	\begin{enumerate}[label=(\roman*)]
		\item\label{tf:stationary} the chain has a unique stationary distribution, and
		\item\label{tf:speed} the average speed of $\froggie_m$, which is defined as the limit
			\[
				s_m\eqdef \lim_{n\to\infty}\frac{\E_{R\sim \bet^n} D_m(F_0,R)}{n},
			\]
			exists and is independent of the initial state~$F_0$, and
		\item\label{tf:constant} for every $\rho\geq 0$, the constant $\gamma_W=\gamma_W(\rho)$ from \Cref{thm:basic} can be expressed as
			\[
				\gamma_W=\rho-\frac{1}{k}\sum_{s_m\leq \rho} (\rho-s_m),
			\]
			and the constant $\tau_W=\tau_W(\rho)$ is nonzero if and only if $\rho=s_m$ for some~$m$.
	\end{enumerate}
\end{theorem}

\paragraph{Special case.}
An interesting special case is when $W$ contains every letter of the alphabet precisely once.
Since the nature of the alphabet is unimportant for us, we may assume that $\bet=[k]$ and
$W=12\dotsb k$. This case admits an elegant closed-form solution.
\begin{theorem}\label{thm:nice}
	Let $\bet=[k]$, $W=12\dotsb k$ and consider the associated frog dynamics. Then the frogs' speeds satisfy $s_m=(k+1)/(k+2-m)(k+1-m)$.

	In particular, for $\rho=1$, we have $\gamma_W=\min_{t\in \Z_+} \frac{k+t^2}{k(t+1)}$, and $\tau_W$ is nonzero precisely when $k$
	is of the form $r^2+r-1$ for some $r\in \Z_+$.
\end{theorem}
Interestingly, the proof of \Cref{thm:nice} does not require computing the stationary distribution of the frog dynamics.
It turns out that $\sum_{m\leq M} s_m$ can be computed
from the simpler chain that is obtained from the frog dynamics by ignoring $\froggie_{M+1},\dotsc,\froggie_k$,
and suppressing the distinction among $\froggie_1,\dotsc,\froggie_M$. This is similar to the arguments in \cite{spitzer}.
The details and the proof of \Cref{thm:nice} are in \Cref{sec:linear}.

However, the stationary distribution of this chain can be described explicitly.
We shall do this by giving the distribution of $\froggie_{m+1}$ conditional on the known
positions of $\froggie_1,\dotsc,\froggie_m$. Recall that a frog arrangement formally is a bijection
$F\colon\{\froggie_1,\dotsc,\froggie_k\}\to\{\lilypad_0,\dotsc,\lilypad_{k-1}\}$.
For brevity, write $\froggie_{\leq m}\eqdef \{\froggie_1,\dotsc,\froggie_m\}$ and $F(\froggie_{\leq m})\eqdef\{F(\froggie_1),\dotsc,F(\froggie_m)\}$.

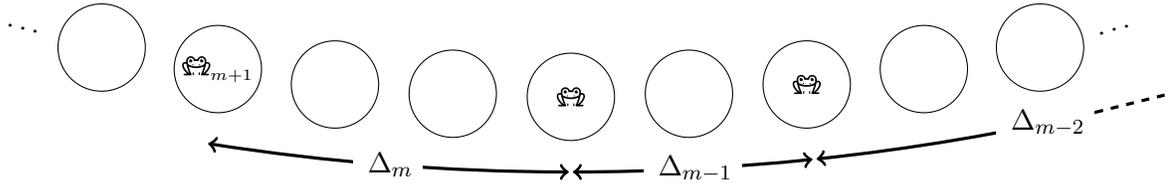
\begin{figure}[ht]
	\centering
        \begin{tikzpicture}
\tikzstyle{lily}=[draw, circle, scale=3.2]
\def\r{30}
\def\rout{31}
\foreach \ang [count=\i from 0] in {-78,-81,...,-102}
  \node[lily] (l\i) at (\ang:\r) {};
\node at (l7) {$\froggie_{\scriptscriptstyle m+1}$};
\node at (l4) {\froggie};
\node at (l2) {\froggie};
\draw [<->, very thick] (-98.95:\rout) arc [start angle=-98.95, end angle=-90.05, radius=\rout] node[pos=0.5,fill=white] {$\,\Delta_m$};
\draw [<->, very thick] (-89.95:\rout) arc [start angle=-89.95, end angle=-84.05, radius=\rout] node[pos=0.5,fill=white] {$\,\Delta_{m-1}$};
\draw [<-, very thick] (-83.95:\rout) arc [start angle=-83.95, end angle=-77, radius=\rout] node[pos=0.82,fill=white] {$\,\Delta_{m-2}$};
\draw [very thick, dashed] (-77:\rout) -- (-75.2:\rout);
\node [rotate=-14]  at (-104:\r) {$\cdots$};
\node [rotate=14]  at (-76:\r) {$\cdots$};
\end{tikzpicture}
        \caption{\label{fig:dyck}Fragment of the lily pad circle in \Cref{thm:margins}. Frogs $\froggie_1$ through $\froggie_m$ are unlabelled.  Note that the gaps are numbered in negative (=clockwise) direction.}
\end{figure}

\begin{theorem}\label{thm:margins}
	Let $F$ be a frog arrangement sampled according to the stationary distribution of the frog dynamics associated with~$W=12\dotsb k$.
	Let $\ell_{m+1}<\ell_m<\ell_{m-1}<\dotsb<\ell_1<\ell_{m+1}+k$, and set $\Delta_i=\ell_{i}-\ell_{i+1}$ for $i\in[m]$. Then
	\[
		\Pr\bigl[F(\froggie_{m+1})=\lilypad_{\ell_{m+1}}\ \big| \ F (\froggie_{\leq m}) = \{\lilypad_{\ell_1}, \dotsc, \lilypad_{\ell_m}\} \bigr]=\frac{1}{\binom{k}{m+1}}\sum_{\substack{a_1,\dotsc,a_m\geq 0\\\sum_{i\leq j} a_i\leq j\\\sum_{i\leq m} a_i=m}}\prod_{i=1}^m\binom{\Delta_i}{a_i}.
	\]
\end{theorem}
The theorem is illustrated in \Cref{fig:dyck}.

\medskip

The formula in \Cref{thm:margins} indicates a curious relationship between the frog dynamics and Dyck paths.
Indeed, consider placing a $+1$ on $\lilypad_{\ell_{m+1}}$ and a $-1$ on each of $\lilypad_{\ell_1},\dots,\lilypad_{\ell_m}$.
Then the sum counts the number of ways to distribute $m$ many $+1$'s onto the $k$ lily pads so that:
\begin{enumerate}
    \item No two $+1$'s occupy the same lily pad, and
    \item All anti-clockwise partial sums starting at $\lilypad_{\ell_{m+1}}$ are strictly positive.
\end{enumerate}
The proof, which is found in \Cref{sec:margins}, exhibits a coupling between the frog dynamics and a particular Markov chain on these arrangements of $\pm 1$'s.

\medskip

The paper is organized as follows.
In \Cref{sec:heights}, we establish the connection between the frog dynamics and the LCS with a periodic word.
We then prove Theorems~\ref{thm:basic}~and~\ref{thm:frog} in \Cref{sec:dynamics}.
\Cref{sec:nice} is devoted to analyzing the case when $W=12\cdots k$: \Cref{sec:linear} contains the proof of \Cref{thm:nice}, and \Cref{sec:margins} contains the proof of \Cref{thm:margins}.
Finally, \Cref{sec:computerjoin} discusses both the new heuristic algorithm for computing the LCS of a pair of random words and the new conjectures which were suggested by our computer simulations.
We conclude the paper with general remarks in \Cref{sec:remarks}.

\medskip

\noindent\textbf{Acknowledgments.} We thank Tomasz Tkocz for discussions at the early stage of this research and for comments on a draft of this paper.
We thank him additionally for the contribution of \Cref{prop:notgaus}.
We thank Alex Tiskin for pointing out the relevance of references~\cite{bundschuh2001} and~\cite{schimd_edit}.
We owe the development of the frog metaphor used in this paper to a conversation with Laure Bukh. The frog symbol
is from Froggy font by Vladimir Nikolic\footnote{Available at \url{https://www.dafont.com/froggy.font}}. The lily pad symbol is based on a drawing
by FrauBieneMaja\footnote{Available at \url{https://pixabay.com/vectors/water-lily-lake-water-pond-blossom-4177686/}}. We thank Zimu Xiang
for pointing several typos, and two anonymous referees for valuable feedback on the earlier versions of the paper.

\section{Heights}\label{sec:heights}
\subsection{Basic definitions and properties}
Let $U$ be any (finite or infinite) word. For a non-negative integer $x$, we denote the $x$-letter prefix by $U_{<x}\eqdef U_0U_1\dotsb U_{x-1}$.
We adopt the convention that $U_{<x}=U$ whenever $x>\len U$.

\begin{defn}\label{def:height}
	A \emph{height} is any function $h\colon \Z\to\Z$ satisfying
	\begin{enumerate}
		\item\label{ht:neg} $h(x)=x$ whenever $x\leq 0$, and
		\item\label{ht:inc} $h(x)-h(x-1)\in \{0,1\}$ for all $x\in \Z$, and
		\item\label{ht:max} $h(x)=h(x-1)$ for all sufficiently large $x$.
	\end{enumerate}
\end{defn}
If $U$ is some fixed word, then for each other word $R$, define $h_R$ by
\begin{alignat*}{2}
	h_R(x)&\eqdef \LCS(R,U_{<x}) &\quad&\text{ for }x\geq 0,\\
	h_R(x)&\eqdef x&\quad&\text{ for }x\leq 0.
\end{alignat*}
It is clear that $h_R$ is a height whenever $R$ is a finite-length word.

For a finite word $R=R_0R_1\cdots R_{\ell-1}$, we will analyze the sequence of heights $h_\varnothing,h_{R_0},h_{R_0R_1},\dots,h_{R}$.
Fix a finite word $R$ and a symbol $a\in\bet$.
Observe that $\LCS(Ra,U_{<x})=\LCS(R,U_{<(x-1)})+1$ if $U_{x-1}=a$, and that $\LCS(Ra,U_{<x})=\max\bigl\{\LCS(R,U_{<x}),\LCS(Ra,U_{<(x-1)})\bigr\}$ if $U_{x-1}\neq a$.
By our convention on the value of $h_R(x)$ if $x\leq 0$ or if $x>\len U$, we have
\begin{equation}\label{eq:evolveh}
	h_{Ra}(x)=\begin{cases}
		h_R(x-1)+1 & \text{if }U_{x-1}=a,\\
		\max\{h_R(x),h_{Ra}(x-1)\} & \text{if }U_{x-1}\neq a.
	\end{cases}
\end{equation}

\paragraph{Heights for periodic words.} It can be shown that every height is of
the form $h_R$ for suitable (possibly infinite) words~$R$ and~$U$ over some alphabet.
However, if $U$ is periodic, then this is reflected in $h_R$ for any word $R$.

For $k\in \mathbb{N}$, define the operator $\delta_k$ by $\delta_k h(x)\eqdef h(x)-h(x-k)$.
Directly from \Cref{def:height}, we observe that if $h$ is any height, then
\begin{itemize}
    \item $\delta_k h(x)=k$ for all $x\leq 0$,
    \item $\delta_k h(x)\in\{0,1,\dots,k\}$ for all $x\in\Z$, and
    \item $\delta_k h(x)=0$ for all $x$ sufficiently large.
\end{itemize}

\begin{defn}
	We say that a height $h$ is a \emph{$k$-height} if $\delta_k h$ is a monotone non-increasing function.
\end{defn}
\begin{theorem}
	If $U$ is a (finite or infinite) $k$-periodic word, then $h_R$ is a $k$-height for every finite word~$R$.
\end{theorem}
\begin{proof}
	We must show $\delta_k h_R(x+k-1)\geq \delta_k h_R(x+k)$ for every~$x\in\Z$. We rewrite this equivalently as
	\begin{equation}\label{eq:kheight}
		h_R(x+k)-h_R(x+k-1)\leq h_R(x)-h_R(x-1).
	\end{equation}
	This inequality holds for $x\leq 0$ because $h_R(x)-h_R(x-1)=1$ in this case.
	Similarly, if $x+k>\len U$, then the inequality holds because $h_R(x+k)-h_R(x+k-1)=0$.
	So, it suffices to prove~\eqref{eq:kheight} only for $1\leq x\leq\len U-k$.
	We do this by induction on~$\len R$, with the base case $\len R=0$ being straightforward.

	Suppose~\eqref{eq:kheight} holds for some word $R$, and we wish
	to establish it for~$Ra$ for some symbol $a\in\bet$.
	Define $\Delta(x)\eqdef h_{Ra}(x)-h_R(x)$.
	From \eqref{eq:evolveh}, we observe that $\Delta(x)\in\{0,1\}$ always.
	\begin{claim}\label{claim:delta}
		If~\eqref{eq:kheight} holds for $R$, then $\Delta(x+k)\geq\Delta(x)$ for all $x\in\Z$ with $x\leq\len U-k$.
	\end{claim}
	\begin{proof}
		We prove the claim by induction on $x$ with the case of $x\leq 0$ being immediate since $\Delta(x)=0$ for these values.

		Suppose that $\Delta(x)=1$; we need to show that $\Delta(x+k)=1$ as well.
		From~\eqref{eq:evolveh}, observe that $\Delta(x)=1$ if and only if either
		\begin{align}
			\label{case:equal} U_{x-1}=a &\quad\text{and}\quad h_R(x)=h_R(x-1),\text{ or}\\
			\label{case:unequal} U_{x-1}\neq a &\quad\text{and}\quad h_{Ra}(x-1)=h_R(x)+1.
		\end{align}

		If~\eqref{case:equal} holds, then since $U$ is $k$-periodic and $x\leq\len U-k$, we have $U_{x+k-1}=a$ as well.
		Furthermore, from~\eqref{eq:kheight}, it follows that $h_R(x+k)=h_R(x+k-1)$ and so $\Delta(x+k)=1$.

		If \eqref{case:unequal} holds, then also $U_{x+k-1}\neq a$.
		Beyond this, $h_{Ra}(x-1)=h_R(x)+1\geq h_R(x-1)+1$ and so $\Delta(x-1)=1$.
		By the induction hypothesis, this implies that $\Delta(x+k-1)=1$ as well.
		Furthermore, $h_R(x-1)=h_{Ra}(x-1)-\Delta(x-1)=h_R(x)$, so~\eqref{eq:kheight} implies that $h_R(x+k)=h_R(x+k-1)$.
		Therefore,
		\begin{align*}
			\Delta(x+k) &= h_{Ra}(x+k)-h_R(x+k)=h_{Ra}(x+k)-h_R(x+k-1)\\
						&\geq h_{Ra}(x+k-1)-h_R(x+k-1)=\Delta(x+k-1)=1.\qedhere
		\end{align*}
	\end{proof}

	Suppose now that we wish to establish \eqref{eq:kheight} for $h_{Ra}$
	in place of $h_R$. The only way the inequality can be violated
	is if
	\begin{align}\label{eq:hraxk}
		h_{Ra}(x+k)&=h_{Ra}(x+k-1)+1,\text{ and}\\
		\label{eq:hrax}
		h_{Ra}(x)&=h_{Ra}(x-1);
	\end{align}
	assume that these hold.
	If $\Delta(x+k-1)=1$, then because of \eqref{eq:evolveh} we must have $h_{Ra}(x+k-1)=h_R(x+k-1)+1\geq h_{Ra}(x+k)$, contradicting \eqref{eq:hraxk}.
	Thus, $\Delta(x+k-1)=0$ and so \Cref{claim:delta} implies that $\Delta(x-1)=0$ as well.
	From here, \eqref{eq:hrax} implies that \[h_R(x-1)=h_{Ra}(x-1)=h_{Ra}(x)=h_R(x)+\Delta(x)\geq h_R(x),\] and so $h_R(x)=h_R(x-1)$ and $\Delta(x)=0$.
	With the aid of the induction hypothesis, \eqref{eq:kheight} and $\Delta(x+k-1)=0$ imply that $h_R(x+k)=h_R(x+k-1)=h_{Ra}(x+k-1)$.
	We deduce that the only way for \eqref{eq:hraxk} to hold is if $U_{x+k-1}=a$, implying that $U_{x-1}=a$ as well.
	Since also $h_R(x)=h_R(x-1)$, \eqref{case:equal} implies that $\Delta(x)=1$; a contradiction.
\end{proof}

\subsection{Ledges and frogs}
From now on we regard $k\in\N$ as fixed; all heights will be derived from the periodic word of period $k$.

For a $k$-height $h$, a \emph{ledge} is an integer $x$ such that $\delta_k h(x+1)=\delta_k h(x)-1$.
Since $\delta_k h$ is a non-increasing integer function which varies between $k$ (for $x\leq 0$) and $0$ (for all sufficiently large $x$), there are precisely $k$ ledges: call them $x_1<x_2<\dots<x_k$.
Equivalently, $x_m$ is the largest integer for which $\delta_k h(x)=k-m+1$.

% How the hell do you force no page break after "Lemma 9." otherwise?
%\pagebreak

\begin{lemma}\ \label{lem:frogsareborn}
	\begin{enumerate}
		\item\label{born:distinct} If $h$ is a $k$-height with ledges $x_1<\dots<x_k$, then $x_1,\dots,x_k$ are non-negative and distinct modulo~$k$.
		\item\label{born:reverse} For any $0\leq x_1<\dots<x_k$ which are distinct modulo $k$,
			\[
				h(x)=x-\sum_{i:\ x_i\leq x}\Bigl\lceil{x-x_i\over k}\Bigr\rceil,
			\]
			is the unique $k$-height with ledges $x_1<\dots<x_k$.
	\end{enumerate}
\end{lemma}
\begin{proof}
	\textit{\Cref{born:distinct}}:
	Let $\Delta(x)\eqdef h(x+1)-h(x)$. Then
	\[
		\delta_k h(x)-\delta_k h(x+1)=\Delta(x-k)-\Delta(x).
	\]
	So, $h$ is a $k$-height if and only if $\Delta(x)\leq \Delta(x-k)$ for all~$x$.
	Since $\Delta$ varies from $1$ (for $x\leq 0$) to $0$ (for large enough $x$),
	in every infinite progression with step $k$ there is a unique $x\geq 0$ such that $\Delta(x-k)-\Delta(x)=1$.
	This shows that $x_1,\dots,x_k$ are non-negative and distinct modulo $k$.\medskip

	\textit{\Cref{born:reverse}}:
	We compute
	\begin{equation}\label{eqn:heightdiff}
		h(x+1)-h(x)=1-\sum_{i:\ x_i\leq x+1}\Bigl\lceil{x+1-x_i\over k}\Bigr\rceil+\sum_{i:\ x_i\leq x}\Bigl\lceil{x-x_i\over k}\Bigr\rceil=1-|\{i:x_i\leq x,\ x_i\equiv x\Mod k\}|.
	\end{equation}
	Since $0\leq x_1<\dots<x_k$ are distinct modulo $k$, we have $h(x+1)-h(x)\in\{0,1\}$.
	From \eqref{eqn:heightdiff}, we see also that $h(x)=x$ for $x\leq 0$ and that $h(x)$ is constant for $x>x_k$.
	Therefore, $h$ is a height.
	Now,
	\[
		\delta_kh(x)=k-\sum_{i:\ x_i\leq x}\Bigl\lceil{x-x_i\over k}\Bigr\rceil+\sum_{i:\ x_i\leq x-k}\Bigl\lceil{x-k-x_i\over k}\Bigr\rceil=k-|\{i:x_i<x\}|,
	\]
	so $\delta_k h$ is non-increasing and has ledges $x_1<\dots<x_k$.

	The uniqueness of $h$ follows from the fact that a $k$-height is determined uniquely by its ledges.
\end{proof}

For a $k$-height $h$ with ledges $x_1<\dots<x_k$, define the function $F_h\colon\{\froggie_1,\dots,\froggie_k\}\to\{\lilypad_0,\dots,\lilypad_{k-1}\}$ by
\[
	F_h(\froggie_m)\eqdef\lilypad_{x_m\bmod k}.
\]
For example, the height of the empty word, $h_\varnothing$, has ledges $x_m=m-1$, and so $F_{h_\varnothing}(\froggie_m)=\lilypad_{m-1}$.
Thanks to \Cref{lem:frogsareborn}, if $h$ is a $k$-height, then $F_h$ is a bijection and is thus a frog arrangement.

For a word $R$, we write $F_R$ in lieu of $F_{h_R}$, e.g.\ $F_\varnothing\eqdef F_{h_\varnothing}$.

\subsection{Evolution of \texorpdfstring{$k$-heights}{k-heights}}
Throughout the preceding discussions of $k$-heights, $U$ was assumed to be an arbitrary $k$-periodic word of some length.
From now on, we will fix $U=W^{(\infty)}$ where $W$ is some fixed word of length $k$. We do not assume that $U$ is irreducible.
Observe that $\LCS(R,W^{(x)})=\LCS(R,U_{<x})=h_R(x)$ for any finite word $R$ and positive integer $x$.

We have seen that a $k$-height is uniquely described by its ledges and that the positions of its ledges induce a frog arrangement.
We next describe how the $k$-height $h_R$ changes as we append symbols to~$R$.

Recall that, for a frog arrangement $F$ and a letter $a\in\bet$, the notation $Fa$ denotes
the frog arrangement resulting from poking lily pads labeled $a$, and then
waiting for all of the frogs to come to a rest. More generally, if $R=R_0\dotsb R_{\ell-1}$,
then $FR$ denotes the frog arrangement obtained by first poking lily pads labeled $R_0$, then those labeled $R_1$, etc.
Recall also that $D_m(F,R)$ denotes the total displacement of $\froggie_m$ during this process.

\begin{theorem}\label{thm:evolution}
	Let $R$ and $T$ be any finite words.
	Suppose that $h_R$ has ledges $x_1<\dots<x_k$ and $h_{RT}$ has ledges $y_1<\dots<y_k$.
	Then
	\begin{align}
		\label{evo:frog} F_{RT} &=F_RT,\text{ and} \\
		\label{evo:speed} y_m &= x_m+D_m(F_R,T).
	\end{align}
\end{theorem}

The proof of this theorem will occupy the rest of the section, and is broken into several steps.
Observe that it suffices to prove \Cref{thm:evolution} only when $T=a$ for some $a\in\bet$ since the full claim then follows by induction on $\len T$.
\medskip

Since each transition in the frog dynamics consists of many individual frog hops, we need to relate these intermediate states to the evolution of the $k$-height.

For a $k$-tuple $X=(x_1,\dotsc,x_k)\in\Z^k$ define the function $h[X]\colon \Z\to\Z$ by
\[
	h[X](x)\eqdef x-\sum_{i\colon x_i\leq x} \Bigl\lceil \frac{x-x_i}{k} \Bigr\rceil.
\]
Thanks to \Cref{lem:frogsareborn}, if $x_1<\dots<x_k$ are non-negative and distinct modulo~$k$, then $h[x_1,\dots,x_k]$ is a $k$-height with ledges $x_1<\dots<x_k$.
We will, however, require this definition even when this condition does not hold.
In any case, observe that
\begin{equation}\label{eq:validincr}
	h[X](x+1)-h[X](x)=1-\abs{\{i: x_i\leq x,\ x_i\equiv x\Mod k\}}.
\end{equation}

For any function $h\colon \Z\to\Z$ and a set $S\subseteq \Z/k\Z$, define
\[
	h^S(x+1)\eqdef
	\begin{cases}
		h(x)+1&\text{if }x\in S\Mod k,\\
		\max\bigl\{h(x+1),h^S(x)\bigr\}&\text{if }x\notin S\Mod k.
	\end{cases}
\]
Recalling that $U=W^{(\infty)}$, a comparison with \eqref{eq:evolveh} shows that if $S=\{i : W_i=a\}$, then
$h_R^S=h_{Ra}$.
\medskip

Of course, when relating $h_R$ and $h_{Ra}$, not every choice of $X,S$ can occur. The appropriate
conditions on the pair $(X,S)$ are captured in the following definition.
\begin{defn}
	Given $X=(x_1,\dotsc,x_k)\in \Z^k$ and $S\subseteq \Z/k\Z$, we say that
	the pair $(X,S)$ is \emph{valid} if it satisfies the following conditions:
	\begin{enumerate}[label=(V\arabic*)]
		\item $0\leq x_1\leq x_2\leq \dotsb\leq x_k$.
		\item Every infinite arithmetic progression of the form $x+k\Z$ contains at most $2$
			elements from among $x_1,\dotsc,x_k$. We say $x$ is \emph{of type $t$} if
			$x+k\Z$ contains exactly $t$ elements from among $x_1,\dotsc,x_k$.
		\item \label{prop:valid} For every $z$ of type $2$, there is $x$ of type $0$ with $x<z$ such that whenever $x<y<z$, then $y\notin S\Mod k$ and $y$ is of type $1$.
		\item If $x\in S\Mod k$, then $x$ is of type $1$.
	\end{enumerate}
\end{defn}

In what follows, for $x\in\Z$, we will abuse notation and simply write $x\in S$ to mean $x\in S\Mod k$.
Similarly, for $x\in\Z$, we will write $S\setminus\{x\}$ to denote the set $\{y\in S:y\not\equiv x\Mod k\}$.

Call a valid pair $(X,S)$ \emph{terminal} if $S=\varnothing$ and every $x\in \Z$ is of type $1$.
Note that if $(X,S)$ is terminal, then $h[X]$ is a $k$-height and $h[X]^S=h[X]$.
Given a non-terminal pair $(X,S)$, let $\ell$ be the least index such that either $x_{\ell}\in S$ or
$x_{\ell}\equiv x_i\Mod k$ for some $i<\ell$. Let $\overline{x}_{\ell}$
be the least integer exceeding $x_{\ell}$ such that $\overline{x}_{\ell}\not\equiv x_i\Mod k$ for all $i<\ell$; observe that $\overline x_\ell-x_\ell\leq k$.
Let $\overline{X}$ be obtained from $X$ by replacing
$x_{\ell}$ with $\overline{x}_{\ell}$, and let $\overline{S}=S\setminus\{x_{\ell},\overline{x}_{\ell}\}$.
Note that the pair $(\overline{X},\overline{S})$ is still valid.

\begin{lemma}\label{lem:dropevol}
	Let $(X,S)$ and $(\overline{X},\overline{S})$ be as above and set $h=h[X]$ and $\overline{h}=h[\overline{X}]$.
	Then $h^S=\overline{h}^{\overline{S}}$.
\end{lemma}
Before we begin the proof \Cref{lem:dropevol}, we first remark on how the lemma implies \Cref{thm:evolution}.

If we start with any valid pair $(X,S)$ and iterate the map $(X,S)\mapsto (\overline{X},\overline{S})$,
we eventually end with a terminal pair. Indeed, each application of the map either decreases $\abs{S}$ (if $x_{\ell}\in S$)
or keeps $\abs{S}$ the same and decreases the number of mod-$k$ residue classes that are of type~$2$ (otherwise).

To each valid pair $(X,S)$ we may associate a function $\F\colon \{\froggie_1,\dotsc,\froggie_k\}\to\{\lilypad_0,\dotsc,\lilypad_{k-1}\}$
defined by $\F(\froggie_m)=\lilypad_{x_m\bmod k}$, and interpret $S$ as the set of those lily pads that contain a single frog,
which is agitated.
In this way, the map $(X,S)\mapsto (\overline{X},\overline{S})$ corresponds to a single, intermediate step in the frog dynamics.
In particular, $\froggie_\ell$ is the nastiest frog which is currently agitated, either from being poked (in the case that $x_\ell\in S$) or being scared off by a nastier frog (in the case that $x_\ell\equiv x_i\Mod k$ for some $i<\ell$).
$\froggie_\ell$ will then leap from $\lilypad_{x_\ell\bmod k}$ to $\lilypad_{\hskip0.08em\relax\overline{x}_\ell\bmod k}$, which is the first available lily pad.
\Cref{lem:dropevol} then implies that $F_{Ra}=F_Ra$ for any word $R$ and symbol $a$, establishing \eqref{evo:frog}.
Furthermore, $D_\ell(F_R,a)=\overline x_\ell-x_\ell$ by definition, which establishes \eqref{evo:speed}.
\medskip

In order to tackle \Cref{lem:dropevol}, we must first establish a few relations between $h[X]$ and $h[X]^S$.
\begin{lemma}\label{lem:validprop}
	If $(X,S)$ is a valid pair and $h=h[X]$, then
	\begin{alignat}{2}
		\label{prop:validi}h(x)&\leq h(z)+1&\qquad&\text{whenever } x\leq z+1,\\
		\label{prop:validiii} h^S(x)&\leq h(z)+1&&\text{whenever }x\leq z+1,\\
		\label{prop:validnew}h^S(x)&\leq h^S(y)&&\text{whenever } x\leq y,\text{ and}\\
		\label{prop:validii} h(x)&\leq h^S(y)&&\text{whenever } x\leq y.
	\end{alignat}
\end{lemma}
\begin{proof}
	Let $T_t$ be the number of elements in the interval $[x,z)$ that are of type $t$. From \eqref{eq:validincr} we deduce
	\[
		h(z)-h(x)\geq T_0-T_2.
	\]
	\ref{prop:valid} implies that the elements of type $0$ and of type $2$ alternate, and so $T_2\leq T_0+1$, implying
	$h(x)\leq h(z)+1$.\medskip

	We next prove \eqref{prop:validiii}. The proof is by induction on $x$. Since $h^S(x)=x$ for $x\leq 0$, the basis
	of induction is clear.
	If $x-1\notin S$, then from the induction hypothesis and \eqref{prop:validi} it follows that $h^S(x)=\max\{h(x),h^S(x-1)\}\leq h(z)+1$.
	If $x-1\in S$, then \ref{prop:valid} implies that in the interval $[x-1,z)$ there are at least
	as many integers of type $0$ as of type~$2$. Hence $h(x-1)\leq h(z)$ and so $h^S(x)=h(x-1)+1\leq h(z)+1$
	in this case as well.\medskip

	To prove \eqref{prop:validnew} it is enough to show that $h^S(x)\leq h^S(x+1)$ holds for all $x$.
	If $x\notin S$, then $h^S(x+1)=\max\{h^S(x),h(x+1)\}\geq h^S(x)$.
	If $x\in S$, then $h^S(x)\stackrel{\eqref{prop:validiii}}{\leq} h(x)+1=h^S(x+1)$.\medskip

	Turning to \eqref{prop:validii}, because of \eqref{prop:validnew} it suffices to establish that $h^S(x)\geq h(x)$. If $x-1\in S$, then $h^S(x)=h(x-1)+1\geq h(x)$ by \eqref{prop:validi}.
	If $x-1\notin S$, then $h^S(x)=\max\{h^S(x-1),h(x)\}\geq h(x)$.
\end{proof}
\def\oh{\overline{h}}% Tired of typing these in
\def\oS{\overline{S}}%
\def\ox{\overline{x}}%
\def\oX{\overline{X}}%
\begin{lemma}\label{lem:evolvejump}
	Let $h$ and $\oh$ be as above. Then
	\begin{alignat}{2}\label{eq:bardiff}
		\oh(x)&=
		\mathrlap{\begin{cases}
				h(x)+1&\text{if }x>x_{\ell}\text{ and }x\in (x_{\ell},\ox_{\ell}]+k\Z,\\
				h(x)&\text{otherwise}.
		\end{cases}}
		\intertext{Furthermore, we have}
		% The 1.1em is experimental to make this align with "if" in "cases" environment
		\stepcounter{equation}
		\label{eq:Jh}\tag{\theequation a}  h(x+1)&=h(x)&\hskip1.09em&\text{ \smash{\vrule depth 10pt height 500pt width 0pt}if }x>x_{\ell}\text{ and }x\in (x_{\ell},\ox_{\ell})+k\Z,\\
		\label{eq:Joh}\tag{\theequation b}\addtocounter{equation}{-1}\refstepcounter{equation}  \oh(x+1)&=\oh(x)&\hskip1.09em&\text{ if }x>x_{\ell}\text{ and }x\in (x_{\ell},\ox_{\ell})+k\Z,\\
		\addtocounter{equation}{-1}
		\label{eq:Js}  h^S(x)&=h(x)+1&&\text{ if }x>x_{\ell}\text{ and }x\in (x_{\ell},\ox_{\ell}]+k\Z.
                \refstepcounter{equation}
	\end{alignat}
\end{lemma}
\begin{proof}
	For \eqref{eq:bardiff}, we observe that $\oh(x)-h(x)=0$ whenever $x\leq x_\ell$.
	If $x>x_\ell$, then since $1\leq\ox_\ell-x_\ell\leq k$,
	\[
		\oh(x)-h(x)=\Bigl\lceil {x-x_\ell\over k}\Bigr\rceil-\Bigl\lceil {x-\ox_\ell\over k}\Bigr\rceil,
	\]
	which is $1$ if $x\in(x_\ell,\ox_\ell]+k\Z$ and $0$ otherwise.\medskip

	We next tackle \eqref{eq:Jh} and \eqref{eq:Joh}. Suppose $x>x_{\ell}$ and $x\in (x_{\ell},\ox_{\ell})+k\Z$. By the minimality of $\ox_{\ell}$, the set $\{x,x-k,x-2k,\dotsc\}$
	contains one of $x_1,\dotsc,x_{\ell-1}$, for otherwise we could have chosen a smaller $\ox_{\ell}$ which is congruent
	to~$x$ modulo~$k$.
	Furthermore, because of \ref{prop:valid}, the set $\{x,x-k,x-2k,\dotsc\}$ contains \emph{precisely one} element of $\{x_1,\dots,x_{\ell-1}\}$, and so $h(x+1)=h(x)$ by \eqref{eq:validincr}.
	Similarly we conclude $\oh(x+1)=\oh(x)$.\medskip

	Since $h^S(x)\leq h(x)+1$ by \eqref{prop:validiii}, to establish \eqref{eq:Js} it suffices to prove that
	$h^S(x)\geq h(x)+1$ for the relevant values of~$x$. We do that by induction on $x$. Suppose first that $x\equiv x_{\ell}+1\Mod k$.
	If $x_{\ell}\in S$, then $h(x)=h(x-1)$ and $h^S(x)=h(x-1)+1$. If $x_{\ell}\notin S$, then $h(x)=h(x-1)-1$
	and $h^S(x)=\max\{h^S(x-1),h(x)\}\geq h^S(x-1)\geq h(x-1)=h(x)+1$. That establishes the base.

	For the induction step, assume that $h^S(x)=h(x)+1$ has already been established and that $x,x+1\in (x_\ell,\ox_\ell]+k\Z$.
	From the minimality of $\ell$, and the fact that $\{x,x-k,x-2k,\dotsc\}$ contains precisely one of $x_1,\dotsc,x_{\ell-1}$, we infer that $x\notin S$.
	Hence, $h^S(x+1)=\max\{h^S(x),h(x+1)\}\geq h^S(x)=h(x)+1=h(x+1)+1$.
\end{proof}

\begin{proof}[Proof of \Cref{lem:dropevol}]
	We will prove that $\oh^{\oS}(x+1)=h^S(x+1)$ by induction on $x$, with the base case of $x<0$ being clear.\smallskip

	\textbf{Case $x\in\oS$:} In this case, $x\in S$ holds as well because $\oS\subseteq S$, so we have $h^S(x+1)=h(x)+1$ and $\oh^{\oS}(x+1)=\oh(x)+1$.
	Hence, $\oh^{\oS}(x+1)=h^S(x+1)$ would follow from $\oh(x)=h(x)$.
	The only way that can fail is if $x>x_\ell$ and $x\in (x_{\ell},\overline{x}_{\ell}]+k\Z$ by \eqref{eq:bardiff}.
	We cannot have $x\in (x_{\ell},\ox_{\ell})+k\Z$ because that would contradict the minimality in the choice of $\ox_\ell$.
	We also cannot have $x\equiv \ox_{\ell}\Mod k$ because that would contradict $x\in \oS$.
	So, either way $\oh^{\oS}(x+1)=h^S(x+1)$ in this case.\medskip

	\textbf{Case $x\notin S$, proof of $\oh^{\oS}(x+1)\geq h^S(x+1)$:} In this case $x\notin \oS$ holds as well because of $\oS\subseteq S$, so $h^S(x+1)=\max\{h^S(x),h(x+1)\}$ and $\oh^{\oS}(x+1)=\max\{\oh^{\oS}(x),\oh(x+1)\}$.
	If $h^S(x+1)=h^S(x)$, then, with the help of induction hypothesis, $h^S(x+1)=h^S(x)=\oh^{\oS}(x)\leq \oh^{\oS}(x+1)$.
	Otherwise, $h^S(x+1)=h(x+1)$, and so $h^S(x+1)=h(x+1)\stackrel{\eqref{eq:bardiff}}{\leq} \oh(x+1)\leq \oh^{\oS}(x+1)$.

	\textbf{Case $x\notin S$, proof of $\oh^{\oS}(x+1)\leq h^S(x+1)$:} If $\oh^{\oS}(x+1)=\oh^{\oS}(x)$, then $\oh^{\oS}(x+1)= h^S(x)$ by the induction hypothesis, so $\oh^{\oS}(x+1)\leq\max\{h^S(x),h(x+1)\}=h^S(x+1)$.

	Thus, suppose that $\oh^{\oS}(x+1)=\oh(x+1)$.
	If $\oh(x+1)=h(x+1)$, we are again done because $h(x+1)\leq\max\{h^S(x),h(x+1)\}=h^S(x+1)$.
	Otherwise, $\oh(x+1)\neq h(x+1)$,
	which is to say $x\in [x_{\ell},\ox_{\ell})+k\Z$ and $x\geq x_{\ell}$.

	If $x\equiv x_{\ell}\Mod k$, then $x_{\ell}$ is
	of type~$2$ with respect to $X$ (for otherwise $x\in S$). So, there is $i<\ell$ such that
	$x_i\equiv x_{\ell}\Mod k$. Since $x_i\leq x_{\ell}\leq x$, it follows that $\oh(x+1)=\oh(x)=h(x)$ by \eqref{eq:validincr} and \eqref{eq:bardiff}
	Hence $\oh^{\oS}(x+1)=h(x)\leq h^S(x+1)$ by \eqref{prop:validii}.

	Suppose $x\in (x_{\ell},\ox_{\ell})+k\Z$ and $x>x_{\ell}$. In this case $\oh(x+1)=h(x+1)+1$, $h(x+1)=h(x)$ and $h^S(x+1)=h(x+1)+1$ by \Cref{lem:evolvejump}.
	Hence $\oh^{\oS}(x+1)=\oh(x+1)=h(x+1)+1=h^S(x+1)$.\medskip

	\textbf{Case $x\in S\setminus \oS$, subcase $x<x_{\ell}$:} In this case the set $\{x,x-k,x-2k,\dotsc\}$ contains neither any of
	$x_1,x_2,\dotsc,x_k$ nor $\ox_{\ell}$. Hence, $h(x+1)=h(x)+1$ and $\oh(x+1)=\oh(x)+1$ implying that
	$\oh^{\oS}(x+1)=\oh(x)+1$ and $h^S(x+1)=h(x)+1$ by \eqref{prop:validii} and \eqref{prop:validiii}.
	Since \eqref{eq:bardiff} tells us that $\oh(x)=h(x)$, we conclude that $\oh^{\oS}(x+1)=h^S(x+1)$ holds in this
	case.\medskip

	\textbf{Case $x\in S\setminus \oS$ and $x\geq x_{\ell}$, subcase $x\equiv x_{\ell}\Mod k$:}
	Since $\{x,x-k,x-2k,\dotsc\}$ contains $x_{\ell}$ and $x_{\ell}\in S$, it follows that $h(x+1)=h(x)$.
	Hence $x^S(x+1)=h(x)+1=h(x+1)+1\stackrel{\eqref{eq:bardiff}}{=}\oh(x+1)$.
	In particular,
	\[
		\oh^{\oS}(x+1)=\max\bigl\{\oh(x+1),\oh^{\oS}(x)\bigr\}\geq \oh(x+1)= h^S(x+1).
	\]
	On the other hand, from the induction hypothesis
	$
		\oh^{\oS}(x)=h^S(x)\stackrel{\eqref{prop:validnew}}{\leq} h^S(x+1),
	$
	and therefore $\oh^{\oS}(x+1)=\max\bigl\{\oh(x+1),\oh^{\oS}(x)\bigr\}\leq h^S(x+1)$.\medskip

	\textbf{Case $x\in S\setminus \oS$ and $x\geq x_{\ell}$, subcase $x\equiv \ox_{\ell}\Mod k$:}
	Since $\{x,x-k,x-2k,\dotsc\}$ contains $\ox_{\ell}$ and $\ox_{\ell}\in S$, it follows from~\eqref{eq:validincr} that $h(x+1)=h(x)$.
	Hence \[h^S(x+1)=h(x)+1=h(x+1)+1\stackrel{\eqref{eq:bardiff}}{=}\oh(x+1)+1.\]
	Therefore, $\oh^{\oS}(x+1)\stackrel{\eqref{prop:validiii}}{\leq} \oh(x+1)+1=h^S(x+1)$.
	On the other hand, $h^S(x+1)= h(x)+1=\oh(x)\stackrel{\eqref{prop:validii}}{\leq} \oh^{\oS}(x+1)$.
\end{proof}

We end this section with the key consequence of \Cref{lem:frogsareborn} and \Cref{thm:evolution}.

\begin{theorem}\label{thm:ledgeloc}
	For any finite word $R$, we have
	\[
		h_R(x)=x-\sum_{i:\ x_i\leq x}\Bigl\lceil{x-x_i\over k}\Bigr\rceil,
	\]
	where $x_i=D_i(F_\varnothing,R)+i-1$.
\end{theorem}
\begin{proof}
	Since the $i$'th ledge of $h_\varnothing$ is $i-1$, the $i$'th ledge of $h_R$ is $D_i(F_\varnothing,R)+i-1$ thanks to \Cref{thm:evolution}.
	The formula for $h_R$ then follows immediately from \Cref{lem:frogsareborn}.
\end{proof}

\section{Frog dynamics}\label{sec:dynamics}
This section is devoted to proving Theorems~\ref{thm:basic}~and~\ref{thm:frog}.

Recall that a word $W\in \bet^k$ is said to be \emph{reducible} if it is of the form $W=U^{(k)}$ where $U\in\bet^\ell$ for some integer $\ell<k$ with $\ell\mid k$; that is $W=UU\cdots U$.
Otherwise $W$ is said to be \emph{irreducible}.\footnote{We warn the reader that ``irreducible words'' are unrelated to ``irreducible Markov chains''.}
Throughout the remainder of the paper, $W$ will be a fixed, irreducible word of length $k$ and we will consider the frog dynamics associated with $W$.

The frog dynamics associated with a word $W$ over an alphabet $\bet$ can be described through a random walk on a directed graph.
Using $\cal F$ to again denote the set of all frog arrangements, define the directed graph $G=G(W,\bet)$ on vertex set $\cal F$ which has a directed edge $F_1\to F_2$ whenever there is some $a\in\bet$ with $F_2=F_1a$.
Note that $G$ may have multi-edges and loops, and that each vertex of $G$ has out-degree $\abs{\bet}$.
The frog dynamics associated with $W$ corresponds precisely to the random walk on~$G$ where each edge is traversed with equal probability.

\subsection{Preliminaries}
We begin with an observation about the frogs' movement.

\begin{proposition}\label{prop:nojump}
	Fix $F\in\cal F$ and $a\in\bet$.
	If $\froggie_m$ hopped in the transition from $F$ to $Fa$, then no frog hopped over $\froggie_m$ in this transition.
\end{proposition}
\begin{proof}
	We will prove that if $\froggie_\ell$ hopped over $\froggie_m$ in the transition from $F$ to $Fa$, then $\froggie_m$ never hopped.

	For $i\in[k]$, define the function $F_i\colon\{\froggie_1,\dots,\froggie_k\}\to\{\lilypad_0,\dots,\lilypad_{k-1}\}$ by
	\begin{equation}\label{eqn:intermediate}
		\F_i(\froggie_j)\eqdef\begin{cases}
			(Fa)(\froggie_j) & \text{if }j\leq i,\\
			F(\froggie_j) & \text{if }j>i.
		\end{cases}
	\end{equation}
	Observe that $\F_0=F$ and $\F_k=Fa$, but that, in general, $\F_i$ is \emph{not} a frog arrangement since multiple frogs may occupy the same lily pad.
	We think of $\F_i$ as the intermediate positions of the frogs after $\froggie_i$ has had the chance to jump.

	Fix any $\ell\in[k]$ such that $\froggie_\ell$ hopped in the transition from $F$ to $Fa$; suppose that $F(\froggie_\ell)=\lilypad_x$ and $(Fa)(\froggie_\ell)=\lilypad_{x+y}$ for some $x\in\{0,\dots,k-1\}$ and $y\in[k]$.
	By definition, in $\F_{\ell-1}$, each of $\lilypad_{x+1},\dots,\lilypad_{x+y-1}$ was occupied by a frog nastier than $\froggie_\ell$.
	This implies that for any $m$ for which $\F_{\ell-1}(\froggie_m)\in\{\lilypad_{x+1},\dots,\lilypad_{x+y-1}\}$, we must have $\F_{\ell-1}(\froggie_m)=(Fa)(\froggie_m)$ since $m<\ell$.
	Thus, we need to show that if $(Fa)(\froggie_m)\in\{\lilypad_{x+1},\dots,\lilypad_{x+y-1}\}$, then $\froggie_m$ did not hop in the transition from $F$ to $Fa$.

	Suppose for the sake of contradiction that one of these frogs did hop; let $t\in\{1,\dots,y-1\}$ be the smallest integer for which $(Fa)(\froggie_m)=\lilypad_{x+t}$ and $\froggie_m$ hopped in the transition from $F$ to $Fa$.
	By the definition of $t$, for all $t'\in\{1,\dots,t-1\}$, the frog at $\lilypad_{x+t'}$ in $F$ never hopped; therefore $\froggie_m$ must have hopped over $\lilypad_{x+1},\dots,\lilypad_{x+t-1}$.
	Furthermore, since $m<\ell$, we have $\F_{m-1}(\froggie_\ell)=F(\froggie_\ell)=\lilypad_x$, and so $\froggie_m$ must have hopped over $\froggie_\ell$ as well; a contradiction since $\froggie_m$ is nastier than $\froggie_\ell$.
\end{proof}

The following proposition will not be used until later; however, due to its similarities with the proposition above, it is convenient to include it and its proof here.

\begin{proposition}\label{prop:hit}
	Fix $F\in\cal F$ and $a\in\bet$.
	If $a$ appears in $W$, then $\sum_{m=1}^k D_m(F,a)=k$.
\end{proposition}
\begin{proof}
	Define $\F_0,\F_1,\dots,\F_k$ as in \eqref{eqn:intermediate}.
	If $\froggie_m$ did not hop in the transition, set $L_m=\varnothing$.
	Otherwise, $\froggie_m$ hopped from $\lilypad_x$ to $\lilypad_{x+y}$ for some $x\in\{0,\dots,k-1\}$ and some $y\in[k]$; in this case, set $L_m=\{x+1,\dots,x+y\}$.
	Observe that $D_m(F,a)=|L_m|$.
	We claim the following:
	\begin{itemize}
		\item For any $m\neq\ell\in[k]$, we have $L_m\cap L_\ell=\varnothing$.

			Suppose not and let $m<\ell$ be such that $L_m\cap L_\ell\neq\varnothing$.
			Suppose that $L_m=\{x+1,\dots,x+s\}$ and $L_\ell=\{y+1,\dots,y+t\}$.
			Since these are two cyclic intervals, we have $L_m\cap L_\ell\neq\varnothing$ if and only if either $x+1\in L_\ell$ or $y+1\in L_m$.
			Since $m<\ell$, we cannot have $y+1\in L_m$ or else $\froggie_m$ hopped over $\froggie_\ell$, which is not possible.
			Thus, $x+1\in L_\ell$.
			Since $x\neq y$, this implies also that $x\in L_\ell$, and so $\froggie_\ell$ hopped over $\lilypad_x$.
			As such, there was some $r<\ell$ with $\F_{\ell-1}(\froggie_r)=\lilypad_x$; since $r\neq m$, $\froggie_r$ must have hopped to this lily pad.
			We conclude that $\froggie_\ell$ thus hopped over a frog who had previously hopped, contradicting \Cref{prop:nojump}.

		\item $\bigcup_{m=1}^k L_m=\Z/k\Z$.

			Fix any $x\in\Z/k\Z$ and let $y\in[k]$ be the smallest integer for which the frog sitting on $\lilypad_{x-y}$ in $F$ hopped in the transition.
			Such a $y$ must exist since there is some lily pad labeled $a$.
			If $m\in[k]$ is such that $F(\froggie_m)=\lilypad_{x-y}$, we claim that $x\in L_m$.
			Indeed, suppose that $Fa(\froggie_m)=\lilypad_{x-y+j}$ so that $L_m=\{x-y+1,\dots,x-y+j\}$.
			In $\F_{m-1}$, $\lilypad_{x-y+j}$ was either empty or contained a frog less nasty than $\froggie_m$.
			In either case, the frog sitting on $\lilypad_{x-y+j}$ in $F$ must have hopped at some point, and so $j\geq y$ by definition.
	\end{itemize}

	Therefore, $\sum_{m=1}^kD_m(F,a)=\sum_{m=1}^k|L_m|=|\Z/k\Z|=k$.
\end{proof}

Moving forward, we will abuse notation slightly and write $\lilypad_i+t=\lilypad_{i+t}$ for $t\in\Z_+$.
For example, $F(\froggie_m)+t$ will denote the lily pad that is $t$ hops ahead from
$\froggie_m$ in the arrangement $F$.

For an arrangement $F\in\cal F$ and $i\neq j\in[k]$, define
\[
	\Delta_F(\froggie_i,\froggie_j)\eqdef \min\{ t\in \Z_+ : F(\froggie_i)+t=F(\froggie_j)\},
\]
that is, the circular distance from $\froggie_i$ to $\froggie_j$ in $F$.
Observe that $\Delta_F(\froggie_i,\froggie_j)=k-\Delta_F(\froggie_j,\froggie_i)$.

Recall that $F_\varnothing$ is the frog arrangement with $F_\varnothing(\froggie_m)=\lilypad_{m-1}$ for each $m\in[k]$; that is, $F_\varnothing$ is the frog arrangement associated with the height of the empty word $h_\varnothing$.
Due to \Cref{thm:ledgeloc}, $F_\varnothing$ plays a special role in our analysis.

\begin{lemma}\label{lem:conn}
	Let $W\in \bet^k$ be irreducible and fix any $F_0\in\cal F$.
	Recursively define $R_t$ to be the label of lily pad $F_t(\froggie_1)$ and $F_{t+1}=F_tR_t$.
	Then, there is an integer $\ell$ for which $F_\ell=F_\varnothing$.
	Furthermore, if $F_0(\froggie_1)=\lilypad_0$, then $k\mid\ell$.
\end{lemma}
\begin{proof}
	By definition, $F_{t+1}(\froggie_1)=F_t(\froggie_1)+1$ for every $t\geq 0$.
	In particular, once we show that $F_\ell=F_\varnothing$ for some $\ell$, it will follow that $k\mid\ell$ if $F_0(\froggie_1)=\lilypad_0$.

	We claim first that there is some $T\in\Z_+$ such that $F_{t+1}(\froggie_i)=F_t(\froggie_i)+1$ for every $i\in[k]$ and $t\geq T$.
	For notational convenience, define $\Delta_t(\froggie_i,\froggie_j)=\Delta_{F_t}(\froggie_i,\froggie_j)$.

	Since $\froggie_1$ always hops, $\froggie_2$ can never jump over $\froggie_1$ (\Cref{prop:nojump}), so $\Delta_t(\froggie_1,\froggie_2)$ is monotonically decreasing.
	Therefore, there is some $T_1\in\Z_+$ for which $\Delta_t(\froggie_1,\froggie_2)$ is constant for all $t\geq T_1$.
	This then implies that $F_{t+1}(\froggie_2)=F_t(\froggie_2)+1$ for every $t\geq T_1$.

	Proceeding by induction on $r$, for some $T_r\in\Z_+$, we know that $F_{t+1}(\froggie_i)=F_t(\froggie_i)+1$ for all $t\geq T_r$ and $i\in[r]$.
	Consider now $\Delta_t(\froggie_r,\froggie_{r+1})$ for $t\geq T_r$.
	Since every frog nastier than $\froggie_{r+1}$ hops one lily pad per beat, $\froggie_{r+1}$ cannot jump over any one of them (\Cref{prop:nojump}).
	Therefore $\Delta_t(\froggie_r,\froggie_{r+1})$ is eventually constant and so there is some $T_{r+1}\in\Z_+$ for which $F_{t+1}(\froggie_{r+1})=F_t(\froggie_{r+1})+1$ for all $t\geq T_{r+1}$ as well.\medskip

	Therefore, there is some $T\in\Z_+$ for which $F_{t+1}(\froggie_i)=F_t(\froggie_i)+1$ for all $t\geq T$ and $i\in[k]$.
	By potentially increasing $T$ by at most $k$, we may suppose that $F_T(\froggie_1)=\lilypad_0$.
	We claim that $F_T=F_\varnothing$, which will establish the claim.

	If not, pick the smallest $r\in\{2,\dots,k\}$ for which $F_T(\froggie_r)=F_T(\froggie_i)+1$ for some $i>r$.
	Since each frog hops by one lily pad for all times $t\geq T$, and $\froggie_r$ is nastier than the frog preceding it, this means that $\froggie_r$ must have been agitated by the monster at each of these times and not by another frog.
	In other words, for each $t\geq 0$, the label of lily pad $F_{T+t}(\froggie_1)=\lilypad_t$ must be the same as the label of lily pad $F_{T+t}(\froggie_r)=\lilypad_{t+r-1}$.
	Since $r\neq 1$, this would imply that $W$ is reducible; a contradiction.
\end{proof}

Denote by $\cal F^*$ the set of frog arrangements $F\in \cal F$ for which there is a word $R$ with $F=F_\varnothing R$, and define the corresponding induced subgraph $G^*\eqdef G[\cal F^*]$.
\Cref{lem:conn} shows that $G^*$ is strongly connected and that every vertex of $G$ has a path to some vertex of $G^*$.
In other words, $\cal F^*$ consists of all the recurrent states of the frog dynamics.

\begin{lemma}\label{lem:aprd}
	If $W\in \bet^k$ is irreducible, then $G^*$ is aperiodic.
\end{lemma}
\begin{proof}
	First, starting at $F_\varnothing$, we clearly see that $F_\varnothing W=F_\varnothing$, so $G^*$ contains a closed walk of length~$k$.

	Now, since $W$ is irreducible, there must be some $i\in\{0,\dots,k-1\}$ for which $W_0\neq W_i$; set $F_1=F_\varnothing W_i$, so that $F_1(\froggie_1)=F_\varnothing(\froggie_1)=\lilypad_0$.
        By \Cref{lem:conn} there is a word $R$ such that $F_1R=F_\varnothing$ of length $rk$ for some $r\in \Z_+$. This
        implies that $G^*$ has a closed walk of length $rk+1$, which is coprime with~$k$.
\end{proof}

These two preceding lemmas tell us that the frog dynamics associated with an irreducible word does indeed have a unique stationary distribution, thus establishing \cref{tf:stationary} of \Cref{thm:frog}.

\begin{theorem}\label{thm:stationary}
	If $W$ is irreducible, then the frog dynamics associated with $W$ has a unique stationary distribution which has support $\cal F^*$.\qed
\end{theorem}

The frog dynamics induces an auxiliary Markov chain on the edges of $G$, which we refer to as the \emph{auxiliary frog dynamics}.
The chain has state-space $\cal E\eqdef\cal F\times\bet$, where a state $(F,a)\in \cal E$ is equally likely to transition to each
of $\abs{\bet}$ many states $(Fa,b)$ with $b\in \bet$.

It is routine to check that if $W$ is irreducible and $\pi$ is the unique stationary distribution of the frog dynamics associated with $W$, then $\hat\pi(F,a)\eqdef{1\over |\bet|}\pi(F)$ is the unique stationary distribution of the auxiliary frog dynamics, which has support $\cal E^*\eqdef\cal F^*\times\bet$.

We now prove \cref{tf:speed} of \Cref{thm:frog} in addition to other preliminary results which will be necessary going forward.
Recall that $D_m(F_0,R)$ is the total displacement of $\froggie_m$ as the word $R$ is applied to $F_0$ and that the \emph{speed} of $\froggie_m$ was defined as
\[
	s_m\eqdef\lim_{n\to\infty}{\E_{R\sim \bet^n}D_m(F_0,R)\over n}.
\]

The following lemma records standard facts about the convergence of an observable of a Markov chain to its mean, as applied to $D_m$ and $s_m$.

\begin{lemma}\label{lem:speeds}
	If $W$ is irreducible and $F_0\in\cal F$ is any initial state, then
	\begin{enumerate}[label=(\roman*)]
		\item\label{spd:exists} $s_m=\E_{(F,a)\sim\hat\pi}D_m(F,a)$, and
		\item\label{spd:Ebound} $\bigl|\E_{R\sim \bet^n} D_m(F_0,R)-s_mn\bigr|\leq O(1)$, and
		\item\label{spd:Pbound} For any fixed $\delta>0$, $\Pr_{R\sim \bet^n}\bigl[\bigl|D_m(F_0,R)-s_mn\bigr|\geq\delta n\bigr]\leq e^{-\Omega(n)}$,
	\end{enumerate}
    where the implicit constants hidden in the big-Oh and big-Omega notations depend on $W$, $F_0$ and $\delta$.
\end{lemma}
\begin{proof}
	Suppose that $R=R_0R_1\dotsc R_{n-1}$; then $D_m(F_0,R)=\sum_{t=0}^{n-1}D_m(F_t,R_t)$ where $F_{t+1}=F_tR_t$.
	Now, since $(F_t,R_t)\in\cal E$, we observe that $D_m(F_t,R_t)$ is simply a function on $\cal E$.

	Since the auxiliary frog dynamics admits a unique stationary distribution and $\cal E$ is finite, the strong law of large numbers for Markov chains from~\cite[Theorem 1.10.2]{norris_markov} implies
	\[
		{1\over n}\sum_{t=0}^{n-1}D_m(F_t,R_t)\to \E_{(F,a)\sim\hat\pi}D_m(F,a),
	\]
	almost surely, which establishes \cref{spd:exists} in a stronger form.

	Let $\hat\pi_t$ denote the distribution of $(F_t,R_t)$.
	Since $\hat\pi$ is the unique stationary distribution of the auxiliary frog dynamics, the inequality $|\hat\pi_t(F,a)-\hat\pi(F,a)|\leq e^{-\Omega(t)}$ for any $(F,a)\in\cal E$ follows from~\mbox{\cite[Corollary~2.8]{aldous_markov}}.
	Therefore, \cref{spd:Ebound} follows by the triangle inequality.

	Finally, \cref{spd:Pbound} follows directly from the Chernoff-type bound in~\cite[Theorem 3.1]{CLLMchernoff}.
\end{proof}

\begin{proposition}\label{prop:speedsep}
	If $W\in \bet^k$ is irreducible, then ${1\over |\bet|}=s_1<s_2<\dots<s_k$.
\end{proposition}
\begin{proof}
	Observe that $\froggie_1$ hops if and only if its lily pad is poked, which happens with probability ${1\over |\bet|}$ independently at each time step; thus, $s_1={1\over |\bet|}$.\medskip

	We show now that $s_m<s_{m+1}$.
	Intuitively, $s_m<s_{m+1}$ if and only if $\froggie_{m+1}$ hops over $\froggie_m$ a positive fraction of time steps in the frog dynamics.
	To make this intuition precise, define $J_m(F_0,R)$ to be the number of times that $\froggie_m$ jumps over $\froggie_{m-1}$ while applying $R$ to $F_0$, where $J_1(F_0,R)$ is defined by imagining a frog $\froggie_0$ who sits permanently between $\lilypad_{k-1}$ and $\lilypad_0$.

	\begin{claim}\label{claim:jump}
		For any $R\in \bet^n$ and an initial state $F_0\in\cal F$,
		\[
			J_m(F_0,R)={1\over k}\biggl(D_m(F_0,R)-D_{m-1}(F_0,R)\biggr)+O(1).
		\]
	\end{claim}
	\begin{proof}
		For any $F\in\cal F$ and $a\in\bet$, we observe that
		\[
			\Delta_{Fa}(\froggie_{m-1},\froggie_m)=\Delta_F(\froggie_{m-1},\froggie_m)+D_m(F,a)-D_{m-1}(F,a)-k\cdot J_m(F,a).
		\]
		Summing this over the trajectory of $F_0$ as we apply symbols of $R$ one-by-one, we obtain
		\begin{align*}
			J_m(F_0,R) &=\biggl\lfloor{1\over k}\biggl(D_m(F_0,R)-D_{m-1}(F_0,R)+\Delta_{F_0}(\froggie_{m-1},\froggie_m)-1\biggr)\biggr\rfloor\\
					   &={1\over k}\biggl(D_m(F_0,R)-D_{m-1}(F_0,R)\biggr)+O(1).\qedhere
		\end{align*}
	\end{proof}
	From \Cref{claim:jump} and the strong law of large numbers for Markov chains, we have
	\[
		{s_{m+1}-s_m\over k}=\lim_{n\to\infty}{\E_{R\sim\bet^n}J_{m+1}(F_\varnothing,R)\over n}=\E_{(F,a)\sim\hat\pi}J_{m+1}(F,a),
	\]
	and so $s_m<s_{m+1}$ if and only if $\E_{(F,a)\sim\hat\pi}J_{m+1}(F,a)>0$.
	Because the support of $\hat\pi$ is precisely $\cal E^*=\cal F^*\times\bet$ and $J_{m+1}$ is a non-negative function, it suffices to find some $F^*\in\cal F^*$ and $a^*\in\bet $ for which $J_{m+1}(F^*,a^*)=1$.

	\medskip

	Starting with $F_0=F_\varnothing$, recursively define $R_t$ to be the label of lily pad $F_t(\froggie_{m+1})$ and $F_{t+1}=F_tR_t$.
	Let $T$ be the smallest integer for which $\froggie_{m+1}$ hops over $\froggie_1$ when applying $R_0\cdots R_{T-1}$ to $F_0$ (a priori, $T$ could be $\infty$).
	Since $\froggie_1$ moves by at most one lily pad on each step and $\froggie_{m+1}$ always hops, we observe that $\Delta_{F_t}(\froggie_{m+1},\froggie_1)$ is weakly decreasing for all $t<T$.

	If $T=\infty$, then there would be some $T_0$ for which $\Delta_{F_t}(\froggie_{m+1},\froggie_1)$ is constant for all $t\geq T_0$, implying that $\froggie_1$ must also hop one lily pad at each time $t\geq T_0$.
	Since $\froggie_1$ hops only when its lily pad is poked, this implies that the labels of lily pads $F_{t+T_0}(\froggie_{m+1})=F_{T_0}(\froggie_{m+1})+t$ and $F_{t+T_0}(\froggie_1)=F_{T_0}(\froggie_1)+t$ must be the same for all $t\geq 0$, contradicting the assumption that $W$ is irreducible.

	Therefore, $T$ is finite, so set $F_0'=F_T$.
	If $\froggie_{m+1}$ hopped over $\froggie_m$ in the transition from $F_{T-1}$ to $F_T$, then we are done; thus, suppose otherwise.
	Since $\froggie_m$ can never hop over $\froggie_{m+1}$, we observe that, in $F_0'$, $\froggie_{m+1}$ now resides between $\froggie_1$ and $\froggie_m$.
	From here, repeatedly poke the lily pad of $\froggie_1$.
	Eventually, thanks to \Cref{lem:conn}, we will arrive back at the arrangement $F_\varnothing$ wherein $\froggie_m$ resides between $\froggie_1$ and $\froggie_{m+1}$.
	Since $\froggie_1$ always hopped when moving from $F_0'$ to $F_\varnothing$, $\froggie_m$ could not have hopped over $\froggie_1$ (\Cref{prop:nojump}); hence $\froggie_{m+1}$ must have hopped over $\froggie_m$ as needed.

	\medskip

	In any case, we have located some $F^*\in\cal F$ and $a^*\in\bet$ for which $J_{m+1}(F^*,a^*)=1$.
	By construction, $F^*\in\cal F^*$, and so we have established the claim.
\end{proof}

\subsection{Limiting distribution of the LCS}
In this section, we compute the mean of $\LCS(R,W^{(\rho n)})$ for a random $R$ and determine its limiting distribution.
To achieve that, we will need to show first that the displacement $D_m$ of $\froggie_m$, once properly shifted and scaled, is distributed asymptotically normally.

For random variables $Y,X_1,X_2,\dots$, write $X_n\tod Y$ if the sequence $X_1,X_2,\dots$ converges to $Y$ in distribution.
That is, $\Pr[X_n\leq t]\to\Pr[Y\leq t]$ for all $t\in\R$ for which $\Pr[Y=t]=0$.
We denote by $\cal N(\mu,\sigma^2)$ the Gaussian distribution with mean $\mu$ and variance $\sigma^2$.

For reasons that will be apparent from \Cref{thm:normnotnorm} below, in addition to normality of $D_m$, we will need to know normality also of $D_1+\dotsb+D_b$.
We prove normality of both $D_m$ and $D_1+\dotsb+D_b$ together.
\begin{lemma}\label{lem:norm}
	Let $W\in \bet^k$ be irreducible, $R\sim \bet^n$ be chosen uniformly at random, $F_0\in\cal F$ be any initial state and $a,b\in[k]$.
	Suppose that either
	\begin{enumerate}
		\item\label{norm:equal} $a=b$, or
		\item\label{norm:lessk} $a=1$ and $b<k$, or
		\item\label{norm:equalk} $a=1$ and $b=k$ and there is some letter of $\bet$ which is absent from $W$.
	\end{enumerate}
	Then there is a $\sigma>0$ for which
	\[
		\sqrt{n}\biggl({\sum_{m=a}^b D_m(F_0,R)\over n}-\sum_{m=a}^bs_m\biggr)\tod\cal N(0,\sigma^2)
	\]
	as $n\to\infty$.
\end{lemma}

The novelty of the above lemma is not that $\sum_{m=a}^b D_m(F_0,R)$ is asymptotically normal; this follows immediately from the central limit theorem for Markov chains.
The novelty stems from the fact that the limiting Gaussian has nonzero variance.
In contrast, if every letter of $\bet$ is present in $W$, then for any $R\in\bet^n$ and $F_0\in\cal F$, we have $\sum_{m=1}^k D_m(F_0,R)=kn$ due to \Cref{prop:hit}.
Thus, in this situation, we have $\sqrt n\bigl({\sum_{m=1}^k D_m(F_0,R)\over n}-\sum_{m=1}^k s_m\bigr)\equiv 0$, which is a degenerate Gaussian.

\begin{proof}
	For any $F\in \cal F$ and any word $S$, put $D(F,S)\eqdef \sum_{m=a}^b D_m(F,S)$. Our aim is to show that $D(F_0,R)$ is asymptotically normal with nonzero variance.
	A general result about Markov chains, which we discuss in \Cref{apx:centlim}, shows that this would follow once we show that there are $F_0,F_1\in\cal F^*$ and words $R,S$ such that $F_0R=F_0$, $F_1S=F_1$, and ${1\over\len R}D(F_0,R)\neq{1\over\len S}D(F_1,S)$.

	We begin by observing that $F_\varnothing W=F_\varnothing$ and that ${1\over\len W}D_m(F_\varnothing,W)={k\over k}=1$ for every $m\in[k]$.
	Hence, ${1\over\len W}D(F_\varnothing,W)=b-a+1$.
	We next exhibit an arrangement $F_0\in\cal F^*$ and a word $R$ such that ${1\over\len R}D(F_0,R)\neq b-a+1$.
	We will need to use a different $F_0$ and $R$ depending on which of the three cases above holds.
	\medskip

	\textbf{Cases~\ref{norm:equal}~and~\ref{norm:lessk}:} Here we have either $a=b$ or $a=1$ and $b<k$; we will use the idea from the proof of \Cref{lem:conn}.
	Call all lily pads labeled $W_0$ \emph{safe}. Starting with $F_\varnothing$, consider repeatedly poking the lily pad containing the nastiest frog that occupies an unsafe lily pad.
	If symbol $W_0$ appears $r$ times in $W$, then, by the argument in the proof of \Cref{lem:conn}, eventually frogs $\froggie_1,\dots,\froggie_r$ occupy all safe lily pads and $\froggie_{r+1},\dots,\froggie_k$ each jump to the very next unsafe lily pad whenever $\froggie_{r+1}$ is poked.
	Thus, if $F_0$ denotes this arrangement and $R_t$ denotes the label of $F_t(\froggie_{r+1})$ with $F_{t+1}=F_tR_t$, then setting $R=R_0R_1\cdots R_{k-r-1}$, we have
	\begin{enumerate}
		\item $F_0=F_{k-r}=F_0R$,
		\item ${1\over k-r}D_m(F_0,R)=0$ for all $m\in[r]$, and
		\item ${1\over k-r}D_m(F_0,R)={k\over k-r}$ for all $m\in\{r+1,\dots,k\}$.
	\end{enumerate}
	If $a=b$, then $D=D_m$, and so ${1\over k-r}D(F_0,R)$ is either $0$ or ${k\over k-r}$; in particular ${1\over k-r}D(F_0,R)\neq 1$.
	If $a=1$ and $b<k$, then
	\[
		{1\over k-r}D(F_0,R)=\begin{cases}
			0 & \text{if }b\leq r,\\
			{k\over k-r}(b-r) & \text{if }b>r.
		\end{cases}
	\]
	In either case, ${1\over k-r}D(F_0,R)\neq b-a+1$ since $a=1$ and $b<k$.
	By construction, $F_0\in\cal F^*$.
	\medskip

	\textbf{Case~\ref{norm:equalk}:} Here we have $a=1$, $b=k$ and there is some letter $y\in\bet$ which is absent from $W$.
	Thus, $F_\varnothing y=F_\varnothing$ and $D_m(F_\varnothing,y)=0$ for all $m\in[k]$, and so ${1\over 1}D(F_\varnothing,y)=0\neq b-a+1$.
\end{proof}

\begin{lemma}\label{lem:approx}
	If $W\in \bet^k$ is irreducible, $\rho$ is a positive real number and $R\in \bet^n$, then
	\[
		\LCS(R,W^{(\rho n)})=\biggl(\rho -{1\over k}\lambda_\rho(R)\biggr)n+O(1),
	\]
	where
	\[
		\lambda_\rho(R)\eqdef\sum_{m=1}^k\max\biggl\{0,\ \rho-{D_m(F_\varnothing,R)\over n}\biggr\}.
	\]
\end{lemma}
\begin{proof}
	A direct application of \Cref{thm:ledgeloc} implies that
	\[
		\LCS(R,W^{(\rho n)})=h_R(\rho n)=\rho n-\sum_{m:\ x_m\leq\rho n}\Bigl\lceil{\rho n-x_m\over k}\Bigr\rceil,
	\]
	where $x_m=D_m(F_\varnothing,R)+m-1$, and so
	\begin{align*}
		\LCS(R,W^{(\rho n)}) &= \rho n-\sum_{m=1}^k\max\biggl\{0,\Bigl\lceil{\rho n-x_m\over k}\Bigr\rceil\biggr\}\\
							 &= \rho n-{1\over k}\sum_{m=1}^k\max\bigl\{0,\rho n-D_m(F_\varnothing,R)\bigr\}+O(1)\\
							 &= \biggl(\rho -{1\over k}\lambda_\rho(R)\biggr)n+O(1).\qedhere
	\end{align*}
\end{proof}

As a next step, we approximate $\lambda_\rho(R)$ by a better-behaved random variable.
Recall that \linebreak $s_1<\dots<s_k$ (\Cref{prop:speedsep}), so let $M$ be the smallest index for which $s_M\geq\rho$ (with $M=k+1$ if $\rho>s_k$) and define the random variable
\[
	\lambda_\rho'(R)\eqdef\sum_{m=1}^{M-1}\biggl(\rho-{D_m(F_\varnothing,R)\over n}\biggr)+\mathbf{1}[s_M=\rho]\cdot\max\biggl\{0,\rho-{D_{M}(F_\varnothing,R)\over n}\biggr\},
\]
where $R\sim\bet^n$.

\begin{lemma}\label{lem:nicervar}
	We have $\E_{R\sim\bet^n}\bigl|\lambda_\rho(R)-\lambda_\rho'(R)\bigr|\leq e^{-\Omega(n)}$ as $n\to\infty$.
\end{lemma}
\begin{proof}
	Using \cref{spd:Pbound} of \Cref{lem:speeds}, we bound
	\begin{align*}
		\Pr[\lambda_\rho(R)\neq\lambda_\rho'(R)] &\leq\sum_{m:\ s_m<\rho}\Pr\biggl[\max\biggl\{0,\rho-{D_m(F_\varnothing,R)\over n}\biggr\}\neq\rho-{D_m(F_\varnothing,R)\over n}\biggr]+\\
		&\qquad\qquad +\sum_{m:\ s_m>\rho}\Pr\biggl[\max\biggl\{0,\rho-{D_m(F_\varnothing,R)\over n}\biggr\}\neq 0\biggr]\\
		&=\sum_{m:\ s_m<\rho}\Pr\bigl[D_m(F_\varnothing,R)>\rho n\bigr]+\sum_{m:\ s_m>\rho}\Pr\bigl[D_m(F_\varnothing,R)<\rho n\bigr]\\
		&\leq\sum_{m:\ s_m\neq\rho}\Pr\bigl[|D_m(F_\varnothing,R)-s_m n|>|\rho-s_m|n\bigr]\leq e^{-\Omega(n)}.
	\end{align*}

	Next, $D_m(F_\varnothing,R)\leq kn$ holds for any $m\in[k]$ and $R\in\bet^n$, which implies the crude inequality $|\lambda_\rho(R)-\lambda_\rho'(R)|\leq 2k(\rho+k)$.
	Therefore,
	\[
		\E\bigl|\lambda_\rho(R)-\lambda_\rho'(R)\bigr|\leq 2k(\rho+k)\Pr\bigl[\lambda_\rho(R)\neq\lambda_\rho'(R)\bigr]\leq e^{-\Omega(n)}.\qedhere
	\]
\end{proof}

\begin{theorem}\label{thm:lambdarho}
	If $W\in \bet^k$ is irreducible, then
	\[
		\E_{R\sim \bet^n}\lambda_\rho(R)=\sum_{s_m\leq\rho}(\rho-s_m)+{\sigma\over\sqrt{2\pi n}}+O(1/n),
	\]
	where
	\[
		\sigma = \begin{cases}
			\lim_{n\to\infty}{1\over n}\sqrt{\Var_{R\sim \bet^n}D_m(F_\varnothing,R)} & \text{if }\rho=s_m,\\
			0 & \text{otherwise.}
		\end{cases}
	\]
	Furthermore, if $\rho\in\{s_1,\dots,s_k\}$, then $\sigma>0$.
\end{theorem}
\begin{proof}
	Combining \Cref{lem:nicervar} with \cref{spd:Ebound} of \Cref{lem:speeds}, and letting $M$ be the smallest index with $s_M\geq\rho$, we compute
	\begin{align*}
		\E\lambda_\rho(R) &= \E\lambda_\rho'(R)+O(e^{-\Omega(n)})\\
						 &= \sum_{m=1}^{M-1}\E\biggl[\rho-{D_m(F_\varnothing,R)\over n}\biggr]+\mathbf 1[s_{M}=\rho]\cdot\E\max\biggl\{0,\rho-{D_{M}(F_\varnothing,R)\over n}\biggr\}+O(e^{-\Omega(n)})\\
						 &= \sum_{m=1}^{M-1}(\rho-s_m)+\mathbf 1[s_{M}=\rho]\cdot\E\max\biggl\{0,\rho-{D_{M}(F_\varnothing,R)\over n}\biggr\}+O(1/n).
	\end{align*}
	If $\rho\notin\{s_1,\dots,s_k\}$, then we are done.

	Otherwise, suppose that $s_M=\rho$.
	Here, \Cref{lem:norm} tells us that $\sqrt{n}\bigl({D_M(F_\varnothing,R)\over n}-\rho\bigr)\tod\cal N(0,\sigma^2)$ where $\sigma>0$.
	Certainly,
	\[
		\sigma^2=\lim_{n\to\infty}{1\over n^2}\Var_{R\sim \bet^n}D_M(F_\varnothing,R).
	\]
	Let $\Phi(t)$ denote the cumulative distribution function of $\cal N(0,\sigma^2)$.
	We can apply the Berry--Esseen-type inequality from~{\cite[Theorem C]{BKloeckner}} to bound
	\[
		\sup_{t\in\R}\biggl|\Pr\biggl[\sqrt n\biggl(\rho-{D_M(F_\varnothing,R)\over n}\biggr)\leq t\biggr]-\Phi(t)\biggr|\leq O(1/\sqrt n).
	\]
	From here, it is clear that
	\[
		\E\max\biggl\{0,\rho-{D_M(F_\varnothing,R)\over n}\biggr\}={1\over\sqrt n}\E_{X\sim\cal N(0,\sigma^2)}\max\{0,X\}+O(1/n)={\sigma\over\sqrt{2\pi n}}+O(1/n),
	\]
	which concludes the proof.
\end{proof}

For the final necessary step, we pin down the asymptotic distribution of the random variable
\[
	\Lambda_\rho(R)\eqdef\sqrt n\bigl(\lambda_\rho(R)-\E\lambda_\rho(R)\bigr),
\]
where $R\sim \bet^n$.

\begin{theorem}\label{thm:normnotnorm}
	Let $W\in \bet^k$ be irreducible, $\rho$ be a positive real number and $R\sim\bet^n$.
	As $n\to\infty$, we have the following:
	\begin{enumerate}[label=(\roman*)]
		\item\label{nnn:avoidspeed} If ${1\over |\bet|}<\rho<s_k$ and $\rho\notin\{s_1,\dots,s_k\}$, then $\Lambda_\rho(R)\tod\cal N(0,\sigma^2)$ for some $\sigma>0$.
		\item\label{nnn:bigmiss} If $\rho>s_k$ and there is some symbol of $\bet$ which is absent from $W$, then $\Lambda_\rho(R)\tod\cal N(0,\sigma^2)$ for some $\sigma>0$.
		\item\label{nnn:toobigtoosmall} If either $\rho>s_k$ and every symbol of $\bet$ appears in $W$ or if $\rho<1/|\bet|$, then $\E|\Lambda_\rho(R)|\to 0$.
		\item\label{nnn:hitspeed} If $\rho\in\{s_1,\dots,s_k\}$, then $\Lambda_\rho(R)$ does \emph{not} converge to a Gaussian distribution.

            In fact, $\Lambda_\rho(R)\tod\max\{G_1,G_2\}-\E\max\{G_1,G_2\}$ where $G_1,G_2$ are centered (possibly degenerate) Gaussian random variables with $\Pr[G_1=G_2]=0$.
	\end{enumerate}
\end{theorem}
\begin{proof}
	Define
	\[
		\Lambda_\rho'(R)\eqdef\sqrt n\bigl(\lambda_\rho'(R)-\E\lambda_\rho'(R)\bigr).
	\]
	\Cref{lem:nicervar} tells us that
	\[
		\E\bigl|\Lambda_\rho(R)-\Lambda_\rho'(R)\bigr|\leq 2\sqrt n\E\bigl|\lambda_\rho(R)-\lambda_\rho'(R)\bigr|\leq e^{-\Omega(n)}.
	\]
	Thus, if $\Lambda_\rho'(R)\tod X$ for some random variable $X$, then also $\Lambda_\rho(R)\tod X$.
	As such, throughout the proof, we will work instead with $\Lambda_\rho'(R)$.
	\medskip

	\textit{\Cref{nnn:avoidspeed}:} Since $s_1<\rho<s_k$ and $\rho\notin\{s_1,\dots,s_k\}$, we have
	\[
		\lambda_\rho'(R)=\sum_{m=1}^M\biggl(\rho-{D_m(F_\varnothing,R)\over n}\biggr)=M\rho-{1\over n}\sum_{m=1}^MD_m(F_\varnothing,R),
	\]
	for some $M\in[k-1]$.
	Since $M$ and $\rho$ are fixed, \Cref{lem:norm} implies that $\Lambda_\rho'(R)\tod\cal N(0,\sigma^2)$ for some $\sigma>0$, and so the same is true of $\Lambda_\rho(R)$.
	\medskip

	\textit{\Cref{nnn:bigmiss}:} Here we have
	\[
		\lambda_\rho'(R)=\sum_{m=1}^k\biggl(\rho-{D_m(F_\varnothing,R)\over n}\biggr)=k\rho-{1\over n}\sum_{m=1}^k D_m(F_\varnothing,R).
	\]
	Since $k$ and $\rho$ are fixed, \Cref{lem:norm} implies that $\Lambda_\rho'(R)\tod\cal N(0,\sigma^2)$for some $\sigma>0$; thus the same is true of $\Lambda_\rho(R)$.
	\medskip

	\textit{\Cref{nnn:toobigtoosmall}:} If $\rho<1/|\bet|=s_1$, then $\lambda_\rho'(R)\equiv 0$.
	If $\rho>s_k$ and every symbol of $\bet$ appears in $W$, then \Cref{prop:hit} tells us that $\lambda_\rho'(R)=k\rho-k$ for any $R\in\bet^n$.
	Since $k$ and $\rho$ are fixed, in either case we have $\Lambda_\rho'(R)\equiv 0$, and so $\E|\Lambda_\rho(R)|\leq e^{-\Omega(n)}\to 0$.
	\medskip

	\textit{\Cref{nnn:hitspeed}:} Suppose that $\rho=s_M$ and set
	\[
		X_n=\sum_{m=1}^{M-1}\biggl(\rho-{D_m(F_\varnothing,R)\over n}\biggr),\qquad Y_n=\rho-{D_{M}(F_\varnothing,R)\over n},
	\]
	so that
	\[
		\lambda_\rho'(R)=X_n+\max\{0,Y_n\}=\max\{X_n,X_n+Y_n\}.
	\]
	By appealing to the multivariate central limit theorem for Markov chains (c.f.~\cite[Section 1.8.1]{geyer_markov}), we find that\footnote{For random variables $Y,X_1,X_2,\ldots\in\R^d$, the statement ``$X_n\tod Y$'' means that $\Pr[X_n\in A]\to\Pr[Y\in A]$ for all Borel sets $A\subseteq\R^d$ with $\Pr[Y\in\partial A]=0$.}
	\[
		\sqrt n\biggl(\begin{bmatrix} X_n \\ Y_n\end{bmatrix}-\begin{bmatrix} \E X_n\\ \E Y_n\end{bmatrix}\biggr)\tod\begin{bmatrix} X\\ Y\end{bmatrix},
	\]
	where
	\[
		\begin{bmatrix} X\\ Y\end{bmatrix}\sim \cal N\biggl(\begin{bmatrix} 0\\ 0\end{bmatrix},\begin{bmatrix} \sigma_X^2 & \sigma_{XY}\\ \sigma_{XY} & \sigma_Y^2\end{bmatrix}\biggr)
	\]
	for some $\sigma_X,\sigma_Y,\sigma_{XY}\in\R$.
	Observe that if $M=1$, then $\sigma_X=0$ since $X_n\equiv 0$ for all $n$.
	However, $\sigma_Y\neq 0$ in any situation thanks to \Cref{lem:norm}.

	Since the map $(x,y)\mapsto\max\{x,x+y\}$ is continuous, we thus have
	\[
		\Lambda_\rho'(R)\tod\max\{X,X+Y\}-\E\max\{X,X+Y\}.
	\]

	Suppose first that $M=k$ and every symbol of $\bet$ is present in $W$; then $X_n+Y_n=ks_k-k$ for all $n$ due to \Cref{prop:hit}.
	This implies that $X+Y\equiv 0$, which is a degenerate Gaussian random variable.
	On the other hand, if $M=k$ and some symbol of $\bet$ is absent from $W$ or $M<k$, then $X+Y\sim\cal N(0,\sigma^2)$ for some $\sigma\neq 0$ thanks to \Cref{lem:norm}.

	In any case, $X,X+Y$ are centered (possibly degenerate) Gaussian random variables.
	Even though $X,X+Y$ are likely dependent, since $\sigma_Y\neq 0$, we have $\Pr[Y=0]=0\implies\Pr[X=X+Y]=0$.
	The proof that that $\max\{X,X+Y\}$ is \emph{not} a Gaussian random variable is left to \Cref{prop:notgaus} in \Cref{apx:notgaus}.
\end{proof}

\subsection{Proof of Theorems~\ref{thm:basic}~and~\ref{thm:frog}}
It remains to put all the pieces together.
We begin by proving \Cref{thm:frog}.
\Cref{tf:stationary} of \Cref{thm:frog} is a consequence of \Cref{thm:stationary}, and \cref{tf:speed} is contained in \Cref{lem:speeds}.
\Cref{tf:constant} is a consequence of \Cref{lem:approx} and \Cref{thm:lambdarho}.

We next consider \Cref{thm:basic}.
Parts \ref{part:piecewiselinear} and \ref{part:slope} follow from \cref{tf:constant} of \Cref{thm:frog}.
\Cref{part:positive} is a direct consequence of \Cref{thm:lambdarho}.

Turning to \cref{part:normality} of \Cref{thm:basic}, \Cref{lem:approx} implies that $\LCS(R,W^{(\rho n)})$ is asymptotically normal with linear variance if and only if $\Lambda_\rho(R)\tod\cal N(0,\sigma^2)$ for some $\sigma>0$.
From \Cref{prop:speedsep} and \cref{tf:constant} of \Cref{thm:frog}, it follows that the slope of $\gamma_W(\rho)$ is nonzero if and only if $\rho<s_k$.
Furthermore, $\tau_W=0$ if and only if $\rho\notin\{s_1,\dots,s_k\}$ due to \Cref{thm:lambdarho}.
Thus, the ``if'' direction of \cref{part:normality} of \Cref{thm:basic} follows from parts \ref{nnn:avoidspeed} and \ref{nnn:bigmiss} of \Cref{thm:normnotnorm}, and the ``only if'' direction follows from parts \ref{nnn:toobigtoosmall} and \ref{nnn:hitspeed} of \Cref{thm:normnotnorm}.

Finally, we need to show that there is indeed an algorithm to compute $\gamma_W(\rho)$ and $\tau_W(\rho)$ in order to establish \cref{part:algorithm} of \Cref{thm:basic}.

Firstly, thanks to \Cref{lem:speeds}, we know that $s_m=\E_{(F,a)\sim\hat\pi}D_m(F,a)$ where $\hat\pi$ is the stationary distribution of the auxiliary frog dynamics.
This stationary distribution can be found by solving a system of linear equations.
In particular, letting $\mathbf P\in\R^{\cal E\times\cal E}$ denote the transition matrix of the auxiliary frog dynamics, then $\hat\pi\in\R^{\cal E}$ is the solution to the linear system $\hat\pi^T(\mathbf P-\mathbf I)=0$ and $\hat\pi^T\mathbf 1=1$ where $\mathbf I\in\R^{\cal E\times\cal E}$ is the identity matrix and $\mathbf 1\in\R^{\cal E}$ is the all-ones vector.
Hence, there is an algorithm to compute each $s_m$, and thus an algorithm to compute $\gamma_W$.

Turning to $\tau_W$, we know that $\tau_W(\rho)=0$ unless $\rho\in\{s_1,\dots,s_k\}$, so suppose that $s_m=\rho$.
Then, according to \Cref{thm:lambdarho}, $\tau_W$ is defined in terms of $\sigma_m\eqdef\lim_{n\to\infty}{1\over n}\sqrt{\Var D_m(F_\varnothing,R)}$, so we need to show that there is an algorithm to compute $\sigma_m$ from $W$.
From $\mathbf P$ and $\hat\pi$, we can compute the fundamental matrix of the chain: $\mathbf Z=(\mathbf I-\mathbf P+\mathbf 1\hat\pi^T)^{-1}$.
Next, let $\mathbf T\in\R^{\cal E\times\cal E}$ be the diagonal matrix with entries $\hat\pi$ and compute the matrix $\mathbf\Gamma=\mathbf T\mathbf Z+(\mathbf T\mathbf Z)^T+\hat\pi\hat\pi^T-\mathbf T$.
Finally,~\cite[Theorem 2.7]{aldous_markov} states that $\sigma_m^2=(D_m-s_m\mathbf 1)^T\mathbf\Gamma (D_m-s_m\mathbf 1)$.
We conclude that there is indeed an algorithm to compute $\tau_W$ from $W$.

\section{Words with distinct symbols}\label{sec:nice}

In this section, we give an explicit expression for the linear term of $\E\LCS(R,W^{(\rho n)})$ in the case where $W$ consists of distinct symbols.
We shall determine also an explicit formula for the stationary distribution in this case. We achieve both of these tasks using the same idea:
We focus only on $\froggie_{\leq m}\eqdef\{\froggie_1,\dots,\froggie_m\}$, suppressing the distinction between $\froggie_1$ through $\froggie_m$.
For the simpler task of computing the linear term, this idea suffices; for more involved task of computing the stationary distribution, we will also separately track
the position of $\froggie_{m+1}$. In either case, we ignore all remaining frogs.

Informally, we may imagine the frog dynamics through the eyes of $\froggie_{m+1}$, to whom all stronger frogs look
equally threatening, and who, at the same time, is oblivious even to the existence of the weaker frogs. From this frog's
point of view, the current state of the frog dynamics can be described by a pair $\bigl(F(\froggie_{m+1}),F(\froggie_{\leq m})\bigr)$,
where $F(\froggie_{\leq m})$ is the $F$-image of the set $\froggie_{\leq m}$.

\subsection{Explicit linear term}\label{sec:linear}
The key observation is that the total speed of the $m$ nastiest frogs can be computed by keeping track of only $F(\froggie_{\leq m})$.
While this observation holds for a general word $W$, we give a proof only in the simpler case when $W$ consists of distinct symbols.

Instead of tracking the full frog arrangement, we shall track the positions of only the first $m$ frogs. As we are not interested
in the position of any individual frog, we record the state of the system into a set $S\in \binom{\{\lilypad_0,\dotsc,\lilypad_{k-1}\}}{m}$.
We call such a set $S$ an \emph{$m$-arrangement}, and continue to refer to elements of $S$ as `frogs', despite not knowing their relative nastiness.

For an $m$-arrangement $S$ and $a\in\bet$, let $Sa$ be the $m$-arrangement obtained by poking the lily pad labeled $a$ (of which there is at most one in this situation).
If this lily pad is occupied, then the resident frog hops one lily pad in the positive direction (anti-clockwise), agitating any frog currently occupying that lily pad.
This continues until all frogs are content once more.

For an $m$-arrangement $S$ and $a\in\bet$, define $H(S,a)$ to be the total number of frogs that hopped in the transition from $S$ to~$Sa$.

\begin{lemma}\label{lem:unstrat}
	If $W$ consists of distinct symbols, then for any $F\in\cal F$, $m\in[k]$ and $a\in\bet$, we have $(Fa)(\froggie_{\leq m})=\bigl(F(\froggie_{\leq m})\bigr)a$ and $\sum_{i=1}^m D_i(F,a)=H(F(\froggie_{\leq m}),a)$.
\end{lemma}
The lemma implies that if $F_0,F_1,\dotsc$ is the frog dynamics, then $F_0(\froggie_{\leq m}),F_1(\froggie_{\leq m}),\dotsc$ is a Markov chain on $m$-arrangements with the transitions described above.
In addition, the total displacement of frogs in the original dynamics can be read off from the behavior of the new, simpler chain.
The $m$-arrangement chain is similar to the PushTASEP, which was studied under the name ``long-range TASEP'' by Spitzer \cite{spitzer} (for
the more general PushASEP model, see Borodin and Ferrari~\cite{borodin_pushasep}). The main difference from the PushTASEP
is that the underlying space of our chain is $\Z/k\Z$, and not~$\Z$.
\begin{proof}
	We first observe that the movement of $\froggie_1,\dots,\froggie_m$ is unaffected by the floundering of the less nasty frogs $\froggie_{m+1},\dots,\froggie_k$. So, since we track only the nastiest $m$ frogs, we may pretend as if $\froggie_{m+1},\dots,\froggie_k$ do not exist.

	Now, consider a slight variation on the frog dynamics where, when $\froggie_i$ is agitated, instead of $\froggie_i$ hopping over any nastier frog, it instead hops onto the very next lily pad.
	If that lily pad is empty, then the frog stops, otherwise there is another frog occupying the lily pad.
	If the current resident is less nasty than $\froggie_i$, then the current resident becomes agitated and will hop on the next step.
	Otherwise, $\froggie_i$ remains agitated and will continue to hop.

	This alternative viewpoint is readily observed to be equivalent to the original frog dynamics.
	If we suppress the distinction among $\froggie_1,\dots,\froggie_m$, then with this alternative viewpoint, when a frog arrives at a currently occupied lily pad, one of the two frogs will hop away at the next step.
	Since this is true regardless of their relative nastiness, we have found that $(Fa)(\froggie_{\leq m})=\bigl(F(\froggie_{\leq m})\bigr)a$ for any $F\in\cal F$ and $a\in\bet$.
	Furthermore, since with the alternative viewpoint, $\sum_{i=1}^m D_i(F,a)$ is simply the total number of hops that took place, the second claim is clear as well.
\end{proof}

We can therefore couple the frog dynamics and the $m$-arrangement chain by simply poking the same lily pad in each chain.
That is to say, starting with the state $F_\varnothing$ in the frog dynamics and the state $\{0,\dots,m-1\}$ in the $m$-arrangement chain, we apply the same random word in both processes.

\begin{theorem}\label{thm:nicespeed}
	If $W$ consists of distinct symbols, then for any $i\in[k]$, $\displaystyle s_i={k(k+1)\over |\bet|(k+2-i)(k+1-i)}$.
\end{theorem}
\begin{proof}
        We prove this in two steps. We first show that the $m$-arrangement chain admits a unique stationary distribution, and that distribution is uniform on ${\{\lilypad_0,\dots,\lilypad_{k-1}\}\choose m}$. We then use that to compute individual frogs' speeds.
		Let $G$ be the digraph with vertex set ${\{\lilypad_0,\dots,\lilypad_{k-1}\}\choose m}$ where $S\to S'$ if there is $a\in\bet$ with $Sa=S'$.

	\paragraph{$G$ is weakly connected.}
        Let $S_{\varnothing}=\{\lilypad_0,\dotsc,\lilypad_{m-1}\}$. Starting with any $S\in {\{\lilypad_0,\dots,\lilypad_{k-1}\}\choose m}$, we may first reach $\{\lilypad_{k-m},\dotsc,\lilypad_{k-1}\}$
        by repeatedly poking the frog that is on the lowest-numbered lily pad. From there we may then reach $S_\varnothing$ by poking the leftmost frog $m$ times.

	\paragraph{$G$ is aperiodic.} When $m<k$, there are always $k-m$ unoccupied lily pads. Hence, at every step there is always
	a positive probability of remaining in the current state, should we poke one of those lily pads. If $m=k$, there is only one state, and so $G$ is trivially aperiodic.

	\paragraph{The stationary distribution of the $m$-arrangement chain is uniform.}
	If $m=k$ then the claim is obvious since there is only one state; thus suppose that $m<k$.
	Observe that for each state $S$, $\deg^{\text{out}}(S)=|\bet|$; we argue that $\deg^{\text{in}}(S)=|\bet|$ as well.
	Firstly, pick any $\lilypad_s\in S$ and let $j\in\{0,\dots,k-1\}$ be the largest integer such that $\{\lilypad_{s-j},\lilypad_{s-j+1},\dots,\lilypad_s\}\subseteq S$; in particular, $\lilypad_{s-j-1}\notin S$ since $m<k$.
	Consider $S'=\bigl(S\setminus\{\lilypad_s\}\bigr)\cup\{\lilypad_{s-j-1}\}$; if $a$ is the label of $\lilypad_{s-j-1}$, then $S'a=S$.
	Observe that these in-edges to $S$ are distinct, so we have found $m$ different in-edges to $S$.
	Furthermore, consider any $a\in\bet$ which is \emph{not} a label of a lily pad in $S$; then $Sa=S$, and so we have located an additional $|\bet|-m$ in-edges to $S$.
	Thus, $\deg^{\text{in}}(S)\geq|\bet|$ for every state $S$, and we conclude that $\deg^{\text{in}}(S)=|\bet|$ by the hand-shaking lemma.

	Since $G$ is a regular weakly connected digraph, it is strongly connected. Indeed, if $G$ were not strongly connected, then there would exist a strongly connected component
	whose total in-degree exceeds the total out-degree, contradicting regularity. Since $G$ is strongly connected, aperiodic and regular, we conclude that the
	$m$-arrangement chain admits a unique stationary distribution, which is uniform on ${\{\lilypad_0,\dots,\lilypad_{k-1}\}\choose m}$.

	\paragraph{Computation of individual speeds.}
	Let $\pi$ denote the stationary distribution of the frog dynamics associated with $W$.
	From the coupling in \Cref{lem:unstrat} and the uniformity of the $m$-arrangement chain, we deduce that
	\begin{equation}\label{eq:sumofspeeds}
		\sum_{i=1}^m s_i=\E_{F\sim\pi,a\sim\bet}\sum_{i=1}^m D_i(F,a)=\E_{S\sim{\{\lilypad_0,\dots,\lilypad_{k-1}\}\choose m},a\sim\bet}H(S,a).
	\end{equation}
	Fix $a\in\bet$ and $S\in{\{\lilypad_0,\dots,\lilypad_{k-1}\} \choose m}$.
	If no lily pad has label $a$, then certainly $H(S,a)=0$.
	Otherwise, suppose that, without loss of generality, $\lilypad_0$ has label $a$.
	In this case, for any $x\in[k]$, $H(S,a)\geq x$ if and only if $\lilypad_0,\lilypad_1,\dots,\lilypad_{x-1}\in S$.
	Invoking cyclic symmetry, we thus have
	\begin{align*}
		\Pr[H(S,a)\geq x] &={m\over |\bet|}\Pr_{S\sim{\{\lilypad_0,\dots,\lilypad_{k-1}\}\choose m}}[\lilypad_1,\dots,\lilypad_{x-1}\in S\ |\ \lilypad_0\in S]\\
						  &={m\over |\bet|}\Pr_{S'\sim{\{\lilypad_1,\dots,\lilypad_{k-1}\}\choose m-1}}[\lilypad_1,\dots,\lilypad_{x-1}\in S']={m\over |\bet|}{{k-x\choose m-x}\over{k-1\choose m-1}},
	\end{align*}
	from which we compute
        \begin{equation}\label{eq:sumofspeedsii}
	\begin{aligned}
		\E_{S\sim{\{\lilypad_0,\dots,\lilypad_{k-1}\}\choose m},a\sim\bet}H(S,a) &=\sum_{x\geq 1}\Pr[H(S,a)\geq x]={m\over |\bet|{k-1\choose m-1}}\sum_{x=1}^k{k-x\choose m-x}\\
																			   &={km\over |\bet|(k+1-m)}.
	\end{aligned}
        \end{equation}

	From \eqref{eq:sumofspeeds} and \eqref{eq:sumofspeedsii} we deduce that
	\[
		s_i={ki\over |\bet|(k+1-i)}-{k(i-1)\over |\bet|(k+1-(i-1))}={k(k+1)\over |\bet|(k+2-i)(k+1-i)}.\qedhere
	\]
\end{proof}

\Cref{thm:nice} follows immediately from \Cref{thm:nicespeed} and \Cref{thm:frog}\ref{tf:constant}
by setting $\bet=[k]$ and $W=12\cdots k$.

\subsection{Stationary distribution}\label{sec:margins}
Here we prove \Cref{thm:margins}, which describes the distribution of $\froggie_{m+1}$ conditional on the known positions of $\froggie_{\leq m}$.
We observe that the cases $m+1=k$ and $m+1=1$ are simple thanks to cyclic symmetry.
Consequently, we fix $m\in[k-2]$ throughout this section.

Similarly to what we did to compute the frogs' speeds, our proof will rely on a coupling of the frog dynamics associated with $W=12\cdots k$ with another easier-to-analyze chain.
That chain, which we denote by $P(m+1,m)$, is significantly more complicated than the $m$-arrangement chain which was used previously. Furthermore, in order to analyze $P(m+1,m)$, we will have to consider a related chain $P(m,m+1)$. We introduce them together.

\paragraph{The Markov chain $P(a,b)$.}

For integers $1\leq a,b\leq k$, the state-space of $P(a,b)$ is ${\Z/k\Z\choose a}\times{\Z/k\Z\choose b}$.
We think of this state-space as recording the positions of two different types of frogs.
For a pair of sets $S=(S^+,S^-)\in{\Z/k\Z\choose a}\times{\Z/k\Z\choose b}$, the first set $S^+$ denotes the position of $a$ many \emph{positive frogs}, whereas $S^-$ denotes the positions of $b$ many \emph{negative frogs}.
Note that these two sets might intersect, and so some lily pads might be occupied by both a positive and a negative frog.

In the coupling that we will construct between $P(m+1,m)$ and the frog dynamics, the negative frogs will correspond to $\froggie_1,\dots,\froggie_m$, yet only one positive frog will correspond to $\froggie_{m+1}$.
\medskip

The chain $P(a,b)$ evolves as follows:

Starting with some arrangement $S\in{\Z/k\Z\choose a}\times{\Z/k\Z\choose b}$, uniformly at random poke one of the $a+b$ signed frogs.\footnote{Note that we poke only one frog, even if it shares a lily pad with another.}
The poked frog becomes agitated and wants to hop away.
\begin{enumerate}
	\item If the poked frog is positive and occupies the same lily pad as a negative frog, no frogs hop, the poked frog calms down and no other frog becomes agitated.
              (This is a special rule for the poked frog, which does not apply to other agitated frogs.)
	\item Otherwise, letting $x$ denote the position of the currently agitated frog, the agitated frog hops from $x$ to $x+1$ and:
		\begin{enumerate}
			\item \label{informalrulea} If there is a frog at $x+1$ with the same sign as the current frog, that frog also becomes agitated and will hop on next step.
			\item \label{informalruleb} Else, if the current frog is negative and $x+1$ contains only a positive frog, that positive frog also becomes agitated and will hop on next step.
			\item \label{informalrulec} Otherwise no new frog becomes agitated.
		\end{enumerate}
              After one of \hyperref[informalrulea]{(a)}, \hyperref[informalruleb]{(b)} or \hyperref[informalrulec]{(c)} happen, the agitated frog that hopped from $x$ to $x+1$ calms down.
\end{enumerate}
This process continues until all frogs are content once more.\medskip

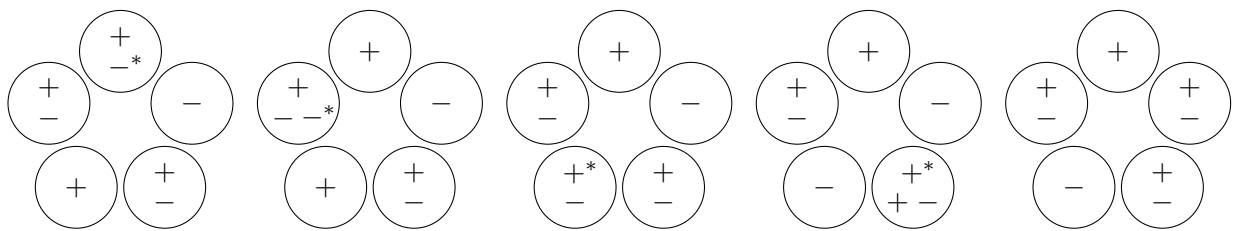
\begin{figure}[ht]
	\begin{center}
		\tikzstyle{lily}=[draw, circle, scale=3]
\def\r{1}
\def\oset{0.2}
\def\sset{0.07}
\begin{tikzpicture}
	\node[lily] (l0) at (90:\r) {};
	\node[lily] (l1) at (162:\r) {};
	\node[lily] (l2) at (234:\r) {};
	\node[lily] (l3) at (-54:\r) {};
	\node[lily] (l4) at (18:\r) {};
	\node at ( $(l0)+(0,\oset)$ ) {$+$};
	\node at ( $(l0)+(\sset,-\oset)$ ) {$-^*$};
	\node at ( $(l1)+(0,\oset)$ ) {$+$};
	\node at ( $(l1)+(0,-\oset)$ ) {$-$};
	\node at (l2) {$+$};
	\node at ( $(l3)+(0,\oset)$ ) {$+$};
	\node at ( $(l3)+(0,-\oset)$ ) {$-$};
	\node at (l4) {$-$};
\end{tikzpicture}\hfil
\begin{tikzpicture}
	\node[lily] (l0) at (90:\r) {};
	\node[lily] (l1) at (162:\r) {};
	\node[lily] (l2) at (234:\r) {};
	\node[lily] (l3) at (-54:\r) {};
	\node[lily] (l4) at (18:\r) {};
	\node at (l0) {$+$};
	\node at ( $(l1)+(0,\oset)$ ) {$+$};
	\node at ( $(l1)+(-\oset,-\oset)$ ) {$-$};
	\node at ( $(l1)+(\oset+\sset,-\oset+0.025)$ ) {$-^*$};
	\node at (l2) {$+$};
	\node at ( $(l3)+(0,\oset)$ ) {$+$};
	\node at ( $(l3)+(0,-\oset)$ ) {$-$};
	\node at (l4) {$-$};
\end{tikzpicture}\hfil
\begin{tikzpicture}
	\node[lily] (l0) at (90:\r) {};
	\node[lily] (l1) at (162:\r) {};
	\node[lily] (l2) at (234:\r) {};
	\node[lily] (l3) at (-54:\r) {};
	\node[lily] (l4) at (18:\r) {};
	\node at ( $(l0)$ ) {$+$};
	\node at ( $(l1)+(0,\oset)$ ) {$+$};
	\node at ( $(l1)+(0,-\oset)$ ) {$-$};
	\node at ( $(l2)+(\sset,\oset)$ ) {$+^*$};
	\node at ( $(l2)+(0,-\oset)$ ) {$-$};
	\node at ( $(l3)+(0,\oset)$ ) {$+$};
	\node at ( $(l3)+(0,-\oset)$ ) {$-$};
	\node at (l4) {$-$};
\end{tikzpicture}\hfil
\begin{tikzpicture}
	\node[lily] (l0) at (90:\r) {};
	\node[lily] (l1) at (162:\r) {};
	\node[lily] (l2) at (234:\r) {};
	\node[lily] (l3) at (-54:\r) {};
	\node[lily] (l4) at (18:\r) {};
	\node at ( $(l0)$ ) {$+$};
	\node at ( $(l1)+(0,\oset)$ ) {$+$};
	\node at ( $(l1)+(0,-\oset)$ ) {$-$};
	\node at (l2) {$-$};
	\node at ( $(l3)+(\sset,\oset)$ ) {$+^*$};
	\node at ( $(l3)-(\oset,\oset)$ ) {$+$};
	\node at ( $(l3)+(\oset,-\oset)$ ) {$-$};
	\node at (l4) {$-$};
\end{tikzpicture}\hfil
\begin{tikzpicture}
	\node[lily] (l0) at (90:\r) {};
	\node[lily] (l1) at (162:\r) {};
	\node[lily] (l2) at (234:\r) {};
	\node[lily] (l3) at (-54:\r) {};
	\node[lily] (l4) at (18:\r) {};
	\node at ( $(l0)$ ) {$+$};
	\node at ( $(l1)+(0,\oset)$ ) {$+$};
	\node at ( $(l1)+(0,-\oset)$ ) {$-$};
	\node at (l2) {$-$};
	\node at ( $(l3)+(0,\oset)$ ) {$+$};
	\node at ( $(l3)+(0,-\oset)$ ) {$-$};
	\node at ( $(l4)+(0,\oset)$ ) {$+$};
	\node at ( $(l4)+(0,-\oset)$ ) {$-$};
\end{tikzpicture}\hfil
	\end{center}
	\caption{\label{fig:signex} An example of one step in the chain $P(4,4)$ with $k=5$. Here, $+$'s indicate positive frogs, $-$'s indicate negative frogs and a ${}^*$ indicates that the frog is agitated. The frogs hop anti-clockwise.}
\end{figure}

Observe that, if the currently agitated frog is positive and hops onto a lily pad occupied only by a negative frog, then this negative frog will \emph{not} become agitated.
In particular, the rules are \emph{not} symmetric between the two signs.
\Cref{fig:signex} displays an example of one step of this process.

In order to analyze $P(a,b)$, we will need to work explicitly with the intermediate steps.
Define $Y^+(a)\eqdef\{1^+,\dots,a^+\}$, $Y^-(b)\eqdef\{1^-,\dots,b^-\}$ and $Y(a,b)\eqdef Y^+(a)\cup Y^-(b)$.
We think of $Y^+$ as the positive frogs and $Y^-$ as the negative frogs.
Let $\Gamma\colon Y(a,b)\to\Z/k\Z$ be any function, which is thought of as an arrangement of signed frogs.
We say that $\Gamma$ is a \emph{valid arrangement} if each lily pad is occupied by at most one frog of each sign.
For a function $\Gamma\colon Y(a,b)\to\Z/k\Z$ and a frog $y\in Y(a,b)$, we say that the pair $(\Gamma,y)$ is a \emph{valid pair} if $\Gamma|_{Y(a,b)\setminus\{y\}}$ is a valid arrangement.

We define three sets:
\begin{itemize}
	\item $\Omega_\beg(a,b)$ is the set consisting of triples $(\Gamma,y,\beg)$ where $\Gamma\colon Y(a,b)\to\Z/k\Z$ is a valid arrangement and $y\in Y(a,b)$.
	\item $\Omega_\nd(a,b)$ is the set consisting of triples $(\Gamma,y,\nd)$ where  $\Gamma\colon Y(a,b)\to\Z/k\Z$ is a valid arrangement and $y\in Y(a,b)$.
	\item $\Omega_\trans(a,b)$ is the set consisting of triples $(\Gamma,y,\trans)$ where $\Gamma\colon Y(a,b)\to\Z/k\Z$ and $y\in Y(a,b)$ with either:
		\begin{itemize}
			\item $\Gamma$ is \emph{not} a valid arrangement, but $(\Gamma,y)$ is a valid pair; or
			\item $\Gamma$ is a valid arrangement, $y\in Y^+(a)$ and there is some $z\in Y^-(b)$ with $\Gamma(z)=\Gamma(y)$.
		\end{itemize}

\end{itemize}
Finally, set $\Omega(a,b)\eqdef\Omega_\beg(a,b)\cup\Omega_\nd(a,b)\cup\Omega_\trans(a,b)$.
Intuitively, for $(\Gamma,y,t)\in\Omega(a,b)\setminus\Omega_\nd(a,b)$, the map $\Gamma$ records the current positions of the signed frogs and $y$ denotes the currently agitated frog.

We turn now to defining a map
\[
	T\colon\Omega(a,b)\setminus\Omega_\nd(a,b)\to\Omega(a,b)\setminus\Omega_\beg(a,b)
\]
which describes the intermediate steps in $P(a,b)$.
Fix $(\Gamma,y,t)\in\Omega(a,b)\setminus\Omega_\nd(a,b)$.
\begin{enumerate}
	\item If $t=\beg$, $y\in Y^+(a)$ and there is some $z\in Y^-(b)$ with $\Gamma(z)=\Gamma(y)$, then $T(\Gamma,y,\beg)=(\Gamma,z,\nd)$.
	\item Otherwise, frog $y$ hops one lily pad forward, which results in $\Gamma'\colon Y(a,b)\to\Z/k\Z$, defined by $\Gamma'(y)=\Gamma(y)+1$ and otherwise agreeing with $\Gamma$.
		\begin{enumerate}
			\item If there is a frog $z\neq y$ of the same sign as $y$ with $\Gamma'(z)=\Gamma'(y)$, then $T(\Gamma,y,t)=(\Gamma',z,\trans)$.
			\item Else, if $y\in Y^-(a)$ and there is $z\in Y^+(a)$ with $\Gamma'(z)=\Gamma'(y)$, then $T(\Gamma,y,t)=(\Gamma',z,\trans)$.
			\item Otherwise, $T(\Gamma,y,t)=(\Gamma',y,\nd)$.
		\end{enumerate}
\end{enumerate}
Observe that $T$ indeed maps $\Omega(a,b)\setminus\Omega_\nd(a,b)$ to $\Omega(a,b)\setminus\Omega_\beg(a,b)$, and hence is well-defined.
Furthermore, for any $(\Gamma,y,\beg)\in\Omega(a,b)$, there is an integer $\ell$ for which $T^{\ell}(\Gamma,y,\beg)\in\Omega_\nd(a,b)$.

We say that $S\in {\Z/k\Z\choose a}\times{\Z/k\Z\choose b}$ and a valid arrangement $\Gamma\colon Y(a,b)\to\Z/k\Z$ are \emph{associated} if $S^+=\Gamma(Y^+(a))$ and $S^-=\Gamma(Y^-(b))$.
Observe that each valid arrangement $\Gamma$ is associated with a unique $S\in {\Z/k\Z\choose a}\times{\Z/k\Z\choose b}$, whereas each $S\in {\Z/k\Z\choose a}\times{\Z/k\Z\choose b}$ is associated with $a!\cdot b!$ valid arrangements.
For an associated $S$ and $\Gamma$, the map $\Gamma$ yields a one-to-one correspondence between $Y(a,b)$ and the signed frogs in $S$.
Hence, we say that $y\in Y(a,b)$ is \emph{associated} with a frog $s$ in $S$ under $\Gamma$.

Now, for any $S\in {\Z/k\Z\choose a}\times{\Z/k\Z\choose b}$ and a signed frog $s$ in $S$, let $S'$ be the result of poking frog $s$.
Select any valid arrangement $\Gamma$ associated with $S$ and let $y\in Y(a,b)$ be the frog associated with $s$ under $\Gamma$.
Since $T$ describes precisely the intermediate steps in $P(a,b)$, if $\ell$ is the integer for which $T^{\ell}(\Gamma,y,\beg)=(\Gamma',y',\nd)$, then $\Gamma'$ and $S'$ are associated.

\paragraph{Coupling $P(m+1,m)$ with the frog dynamics.}
We consider a variant on $P(m+1,m)$, which is slowed down just enough in order to couple it with the frog dynamics.
The chain $\hat P(m+1,m)$ has the same state-space as $P(m+1,m)$ but evolves according to:
\begin{enumerate}
	\item With probability $\max\bigl\{0,1-{2m+1\over k}\bigr\}$, do nothing.
	\item Otherwise, follow the same process as $P(m+1,m)$.
\end{enumerate}
Since $\hat P(m+1,m)$ is simply a (potentially) lazy version of $P(m+1,m)$, any stationary distribution of $P(m+1,m)$ is also a stationary distribution of $\hat P(m+1,m)$.

Consider any state $S=(S^+,S^-)\in{\Z/k\Z\choose m+1}\times{\Z/k\Z\choose m}$.
For $x\in\Z/k\Z$, consider the partial sums
\begin{equation}\label{eqn:partial}
	\sum_{i=x}^{x+j}\bigl(\mathbf{1}[i\in S^+]-\mathbf{1}[i\in S^-]\bigr)\qquad\text{for $j\in\{0,\dots,k-1\}$}.
\end{equation}
It is a well-known fact in the study of Dyck paths (which according to \cite[p.~373]{concretemath} is originally due to Raney~\cite{raney}),
that for any such $S$ there is a \emph{unique} $x\in\Z/k\Z$ for which \eqref{eqn:partial} is strictly positive for all~$j$.
Note that such an $x$ must satisfy $x\in S^+\setminus S^-$.
We call the positive frog sitting at position $x$ \emph{the optimistic frog}.
From the preceding observation about $x$ we know that the optimistic frog does not share its lily pad with another frog.

For a state $S=(S^+,S^-)\in{\Z/k\Z\choose m+1}\times{\Z/k\Z\choose m}$, define two functions $f^-(S)\eqdef\{\lilypad_i:i\in S^-\}$ and $f^+(S)\eqdef\lilypad_x$ where $x$ is the position of the optimistic frog.
Define also $f(S)\eqdef\bigl(f^+(S),f^-(S)\bigr)$.

We say that a state $S\in{\Z/k\Z\choose m+1}\times{\Z/k\Z\choose m}$ and a frog arrangement $F\in \cal F$ are \emph{compatible} if
$f(S)=\nobreak\bigl(F(\froggie_{m+1}),F(\froggie_{\leq m})\bigr)$. Note that many frog arrangements are compatible with a given $S\in {\Z/k\Z\choose m+1}\times{\Z/k\Z\choose m}$, and that a given frog arrangement is compatible with many states in $P(m+1,m)$.

Fix a state $S\in{\Z/k\Z\choose m+1}\times{\Z/k\Z\choose m}$ and let $F\in\cal F$ be any frog arrangement compatible with $S$.
Taking one step in the chain $\hat P(m+1,m)$ results in a new state $S'$.
We couple this action with the frog dynamics in the following way:
\begin{enumerate}
	\item If a negative frog or the optimistic frog was poked in $S$, poke the corresponding lily pad in the frog arrangement.
	\item If a non-optimistic positive frog was poked in $S$ or no frog was poked, uniformly at random select a letter $a\in\bet$ which is not a label of a lily pad in $F(\froggie_{\leq m+1})$ and poke any lily pad labeled $a$.
\end{enumerate}
This will result in a new frog arrangement~$F'$.

We observe first that the above coupling preserves the transition probabilities in the frog dynamics associated with $W=12\cdots k$ and either $\bet=[k]$ if $2m+1\leq k$, or $\bet=[2m+1]$ otherwise.\footnote{Observe that if $|\bet|>k$, then the only difference with the frog dynamics where $\bet=[k]$ is that at each step, the probability that no lily pad is poked is $1-{k\over |\bet|}$. As such, the stationary distribution of the frog dynamics over the larger alphabet is the same as for $\bet=[k]$.}
Indeed, if $2m+1\leq k$, then each lily pad in the frog dynamics is poked with probability $1/k$, and if $2m+1>k$, then each lily pad in the frog dynamics is poked with probability $1/(2m+1)$.

We now verify that the above is indeed a indeed a coupling.

\begin{theorem}
If $S$ and $F$ are compatible, then so are $S'$ and $F'$, i.e.\ $f(S')=\bigl(F'(\froggie_{m+1}),F'(\froggie_{\leq m})\bigr)$.
\end{theorem}
\begin{proof}
	If none of the $2m+1$ frogs in $S$ were poked, then $S'=S$.
	In the coupling, this corresponds to either poking no lily pad or poking a lily pad containing a frog less nasty than $\froggie_{m+1}$; thus $F'(\froggie_i)=F(\froggie_i)$ for all $i\in[m+1]$, as needed.

	Next suppose that some frog in $S$ was poked.
	Since negative frogs move unabated by positive frogs, the movement of the set of negative frogs is identical to the movement of $\{\froggie_1,\dots,\froggie_m\}$ in the $m$-arrangement chain from the previous section.
	Since one of $\froggie_1,\dots,\froggie_m$ moves if and only if one of the negative frogs in $S$ is poked, thanks to \Cref{lem:unstrat}, we know that $f^-(S')=F'(\froggie_{\leq m})$.

	This being the case, it remains to verify only that $f^+(S')=F'(\froggie_{m+1})$.

	Let $\Gamma_0\colon Y(m+1,m)\to\Z/k\Z$ be any valid arrangement associated with $S$ and let $y_0\in Y(m+1,m)$ be the frog associated with the poked frog under $\Gamma_0$.
	Let
	\[
		(\Gamma_0,y_0,\beg),(\Gamma_1,y_1,\trans),\dots,(\Gamma_{\ell-1},y_{\ell-1},\trans),(\Gamma_\ell,y_\ell,\nd)
	\]
	be the trajectory of $(\Gamma_0,y_0,\beg)$ under the map $T$, so that $\Gamma_\ell$ is associated with $S'$.
	Analogously to \eqref{eqn:partial}, for $i\in\{0,\dots,\ell\}$ and $x,x'\in\Z/k\Z$, define
	\[
		L_i[x,x']\eqdef\sum_{r=x}^{x'}\bigl(|\Gamma_i^{-1}(r)\cap Y^+(m+1)|-|\Gamma_i^{-1}(r)\cap Y^-(m)|\bigr).
	\]
    The $L_i[x,x']$ counts the number of frogs in the interval $[x,x']$ weighted by their signs.

    Although $\Gamma_i$ may not be a valid arrangement, by the same reasoning that was applied to \eqref{eqn:partial}, there is still a unique $x_i\in\Z/k\Z$ for which $L_i[x_i,x_i+j]>0$ for all $j\in\{0,\dots,k-1\}$.
	Of course, $x_0$ is the position of the optimistic frog in $S$ and $x_\ell$ is the position of the optimistic frog in $S'$.

	We verify that $F'(\froggie_{m+1})=\lilypad_{x_\ell}$ through the following sequence of claims.
	\begin{claim}\label{claim:nomove}
        For $i\in[\ell]$, if $x_{i-1}\notin\{\Gamma_{i-1}(y_{i-1}),\Gamma_i(y_i)\}$, then $x_i=x_{i-1}$.
	\end{claim}
	\begin{proof}
        Let $\bar x\eqdef\Gamma_{i-1}(y_{i-1})$.
		Suppose first that $y_{i-1}\in Y^-(m)$; then $\Gamma_i(y_{i-1})=\bar x+1=\Gamma_i(y_i)$.
		Since $\Gamma_i(y_i)\neq x_{i-1}$, it follows that $L_i[x_{i-1},x']=L_{i-1}[x_{i-1},x']$ unless $x'=\bar x$ holds, in which
                case $L_i[x_{i-1},x']=L_{i-1}[x_{i-1},x']+1$. In particular, $L_i[x_{i-1},x']\geq L_{i-1}[x_{i-1},x']>0$ for all $x'$, and so $x_i=x_{i-1}$.

		Suppose next that $y_{i-1}\in Y^+(m+1)$; here we have three cases.
		\begin{enumerate}
			\item $i=1$ and there is some $z\in Y^-(m)$ with $\Gamma_0(z)=\Gamma_0(y_0)$:
				In this case, no frogs hop, implying that $\Gamma_1=\Gamma_0$, so $x_1=x_0$.
			\item $i=1$ and there is no $z\in Y^-(m)$ with $\Gamma_0(z)=\Gamma_0(y_0)$:
                                Since the only difference between $\Gamma_0$ and $\Gamma_1$ is that $y_0$ has hopped forward from $\Gamma_0(y_0)$, we have $L_1[x_0,x']=L_0[x_0,x']$ unless $x'=\bar x$.
                                Furthermore, since lily pad $\bar x$ contains no negative frogs, we have $L_1[x_0,\bar x]\geq L_1[x_0,\bar x-1]>0$.
                                So, $L_1[x_0,x']>0$ either way.
			\item $i\geq 2$:
				Again, $y_{i-1}$ hops to the very next lily pad, so $\Gamma_i(y_{i-1})=\bar x+1=\Gamma_i(y_i)$ and otherwise $\Gamma_i$ and $\Gamma_{i-1}$ agree.
				Since $L_i[x_{i-1},x']\geq L_{i-1}[x_{i-1},x']$ unless $x'=\bar x$,
                                the only way for $L_i[x_{i-1},x']\leq 0$ to happen is if $x'=\bar x$.
				Here we must break into cases depending on the sign of $y_{i-2}$.
				\begin{enumerate}
					\item If $y_{i-2}\in Y^-(m)$, then we must have had
                                          \begin{align*}
                                          L_{i-1}[x_{i-1},\bar x-1]&=L_{i-1}[x_{i-1},\bar x]-L_{i-1}[\bar x,\bar x]\\&=L_{i}[x_{i-1},\bar x]\\&=0;
                                          \end{align*}
                                          contradicting the definition of $x_{i-1}$.
					\item Otherwise, $y_{i-2}\in Y^+(m+1)$, which implies that
                                          $L_i[\bar x,\bar x]\geq 0$ and so \[\Gamma_{i-1}[x_{i-1},\bar x]=\Gamma_{i-1}[x_{i-1},\bar x-1]+L_i[\bar x,\bar x]>0.\qedhere\]
				\end{enumerate}
		\end{enumerate}
	\end{proof}
	\begin{claim}\label{claim:negpush}
        Suppose that $x_0\in\bigl\{\Gamma_i(y_i):i\in\{0,\dots,\ell\}\bigr\}$ and let $i$ be the smallest index for which $\Gamma_i(y_i)=x_0$.
		If $i\geq 1$, then $y_{i-1}\in Y^-(m)$.
	\end{claim}
	\begin{proof}
		Suppose not, so $y_{i-1}\in Y^+(m+1)$ and let $r$ be the smallest index for which $y_r\in Y^+(m+1)$.
		Since positive frogs cannot agitate negative frogs, we observe that $y_r,y_{r+1},\dots,y_{i-1}\in Y^+(m+1)$.
                Note that frogs $y_r,y_{r+1},\dotsc,y_{i-1},y_i$ occupy consecutive positions under $\Gamma_r$, and so, letting $\bar x=\Gamma_r(y_r)$, we have
\begin{equation}\label{eq:negpushaux}
               L_r[\bar x,x']=L_r[y_r,x_0-1]+L_i[x_0,x']\text{ for }x'>x_0.
\end{equation}
		We proceed with three cases.
		\begin{enumerate}
			\item $r=0$ and there is some $z\in Y^-(m)$ with $\Gamma_0(z)=\Gamma_0(y_0)$:
				Here, no frogs hop, implying that having $\Gamma_i(y_i)=x_0$ is impossible since $i\geq 1$.
			\item $r=0$ and there was no $z\in Y^-(m)$ with $\Gamma_0(z)=\Gamma_0(y_0)$:
				In particular, $L_0[\bar x,x']>0$ for $\bar x\leq x'<x_0$. In view of \eqref{eq:negpushaux}, this implies that $L_0[\bar x,x']>0$ for all $x'$, contradicting the definition of $i$.
			\item $r\geq 1$:
				Here, by the definition of $r$, $y_r$ must have been agitated by a negative frog.
				Therefore, there was no negative frog occupying the same lily pad as $y_r$ in $\Gamma_{r-1}$.
				Similarly to the previous case, this implies that $L_{r-1}[\bar x,x']>0$ for $\bar x\leq x'< x_0$, which in turn implies
                                that $L_{r-1}[\bar x,x']>0$ for all $x'$. Hence, $x_{r-1}=\bar x$.
				However, thanks to \Cref{claim:nomove}, we know that $x_{r-1}=x_0$; again contradicting the definition of $i$.\qedhere
		\end{enumerate}
	\end{proof}

	With the help of the above claims, we can now deduce that $F'(\froggie_{m+1})=\lilypad_{x_\ell}$.

    Suppose first that the optimistic frog was never agitated, so $x_0\notin\bigl\{\Gamma_i(y_i):i\in\{0,\dots,\ell\}\bigr\}$.
    Then, thanks to \Cref{claim:nomove}, we have $x_0=x_1=\dots=x_\ell$.
    Furthermore, $F'(\froggie_{m+1})=F(\froggie_{m+1})$ and so the claim follows.

    On the other hand, suppose that the optimistic frog was agitated at some point, which is to say $x_0\in\bigl\{\Gamma_i(y_i):i\in\{0,\dots,\ell\}\bigr\}$.
    Let $i$ be the smallest index for which $\Gamma_i(y_i)=x_0$ and consider the largest $j\geq 1$ for which $\{x_0+1,\dots,x_0+j-1\}\subseteq\Gamma_i(Y^-(m))$.
    In the frog dynamics, we would have $F'(\froggie_{m+1})=\lilypad_{x_0+j}$, so we must show that $x_\ell=x_0+j$.

    If $i=0$, then the optimistic frog was poked.
    Otherwise $i\geq 1$: appealing to \Cref{claim:negpush} and using the fact that a positive frog cannot agitate a negative frog, we know that $y_0,\dots,y_{i-1}\in Y^-(m)$, and \Cref{claim:nomove} implies that $x_{i-1}=x_0$.
    So, both in the case $i=0$ and in the case $i\geq 1$, no other positive frog was agitated prior to the optimistic frog, and hence $\Gamma_0(Y^+(m+1))=\Gamma_i(Y^+(m+1))$; in particular, we must have $\{x_0+1,\dots,x_0+j-1\}\subseteq\Gamma_i(Y^+(m+1))$ as well.
    Using these observations, we see that $x_{i+j}=x_0+j$ since, by assumption, $x_0+j\notin\Gamma_i(Y^-(m))=\Gamma_{i+j}(Y^-(m))$.
    If $i+j=\ell$, then we are done.
    Otherwise, position $x_{i+j}$ is occupied by two positive frogs in $\Gamma_{i+j}$.
    In this case, it is easy to observe that $x_{i+j}=x_{i+j+1}$, and we then conclude that $x_{i+j}=x_{i+j+1}=\dots=x_\ell$ by a final appeal to \Cref{claim:nomove}.
\end{proof}

\paragraph{A stationary distribution of $P(a,b)$.}
We show now that for any $1\leq a,b\leq k$, the chain $P(a,b)$ admits \emph{a} stationary distribution which is uniform on ${\Z/k\Z\choose a}\times{\Z/k\Z\choose b}$.
Note that we claim only that \emph{a} stationary distribution of this form exists, not that it is unique.\footnote{In the special case of $P(m+1,m)$ that we care about, one can show that the stationary distribution is indeed unique, but this fact is unnecessary for our arguments.}

Let $G$ be the digraph whose vertices are the states of $P(a,b)$ where $S\to S'$ if poking some frog in $S$ results in $S'$.
In other words, $P(a,b)$ is the random walk on $G$ where each edge is traversed with equal probability.

Decompose $G$ into its weakly connected components $G=G_1\cup\dots\cup G_\ell$.

\begin{lemma}\label{lem:aperiod}
	Each $G_i$ is aperiodic.
\end{lemma}
\begin{proof}
	We claim that there is some $S\in V(G_i)$ with $S^+\cap S^-\neq\varnothing$.
	Indeed, pick some $S_0\in V(G_i)$ with $S_0^+\cap S_0^-=\varnothing$. By starting from $S_0$ and repeatedly poking a fixed positive frog, we will eventually drive that frog to a lily pad occupied by a negative frog. The resulting state $S$ satisfies $S^+\cap S^-\neq\varnothing$ and $S\in V(G_i)$.

    Thus, pick a state $S\in V(G_i)$ and a lily pad $x$ such that $x\in S^+\cap S^-$.
	Poking the positive frog in position $x$ leaves the state unchanged, so $G_i$ contains a closed walk of length $1$.
\end{proof}

\begin{lemma}\label{lem:regular}
	For every vertex $S$ of $G$, $\deg^{\text{out}}(S)=\deg^{\text{in}}(S)=a+b$.
\end{lemma}
\begin{proof}
	We already know that $\deg^{\text{out}}(S)=a+b$.
	In order to show that $\deg^{\text{in}}(S)=a+b$, we establish a much stronger property of $P(a,b)$.

	We define a ``reversal map'' $R\colon\Omega(a,b)\to\Omega(b,a)$.
	Intuitively, the map $R$ will switch the signs of the frogs and reverse the direction of the ring of lily pads.

	Fix $(\Gamma,y,t)\in\Omega(a,b)$; we will define $R(\Gamma,y,t)$. First, define the function $\Gamma'\colon Y(b,a)\to\Z/k\Z$ by $\Gamma'(i^\pm)=k-\Gamma(i^\mp)\bmod k$.
	\begin{itemize}
		\item If $t\in\{\beg,\nd\}$, let $t'$ be such that $\{t,t'\}=\{\beg,\nd\}$. Then $R(\Gamma,i^\pm,t)=(\Gamma',i^\mp,t')$.
		\item Otherwise, $t=\trans$ and so there is some $z\neq y$ with $\Gamma(z)=\Gamma(y)$ where either $z$ and $y$ have the same sign or $y\in Y^+(a)$ and $z\in Y^-(b)$.
			\begin{itemize}
				\item If $y=i^\pm$ and $z=j^\pm$, then $R(\Gamma,i^\pm,\trans)=(\Gamma',j^\mp,\trans)$.
				\item If $y=i^+$ and $z=j^-$, then $R(\Gamma,i^+,\trans)=(\Gamma',j^+,\trans)$.
			\end{itemize}
	\end{itemize}
	Observe that $R$ is an involution.

	We prove now the key time-reversal property of $P(a,b)$: the map $R$ reverses the flow of time, making $T$ its own inverse.
	\begin{claim}\label{claim:reverse}
		The map $RTRT\colon\Omega(a,b)\setminus\Omega_\nd(a,b)\to\Omega(a,b)\setminus\Omega_\nd(a,b)$ is the identity map.
	\end{claim}
	\begin{proof}
		We observe that $R$ maps $\Omega_\trans(a,b)$ to $\Omega_\trans(b,a)$ and swaps $\Omega_\beg(a,b)$ and $\Omega_\nd(b,a)$, so the map $RTRT$ is well-defined.

		For $(\Gamma,y,t)\in\Omega(a,b)\setminus\Omega_\nd(a,b)$, observe that if $\Gamma(y)=x$, then $T$ affects only frogs on lily pads $x$ and $x+1$; hence, we need keep track of only these two lily pads.
		This means that we may prove the claim by checking every possible arrangement of positive and negative frogs on two consecutive lily pads.
		That amounts to straightforward, but slightly tedious case-checking, which we defer to \Cref{apx:reverse}.
	\end{proof}
	A corollary of the above claim is that for any $\ell\geq 1$, $RT^{\ell}RT^{\ell}=\text{id}$ as well.
	Indeed, using the fact that $R$ is an involution and proceeding by induction on $\ell$ we have
	\[
		RT^{\ell}RT^{\ell}=RTR\bigl( RT^{\ell-1}RT^{\ell-1}\bigr)T=RTRT=\text{id}.
	\]

	Let $S\in {\Z/k\Z\choose a}\times{\Z/k\Z\choose b}$ be arbitrary, and let $\Gamma$ be any valid arrangement associated with $S$.
	Suppose $S'\in {\Z/k\Z\choose a}\times{\Z/k\Z\choose b}$ and $s'$ is a signed frog so that poking $s'$ in $S'$ results in~$S$.
	Then there is a unique $\Gamma'$ associated with $S'$ and $y'$ associated with $s'$ under $\Gamma'$ for which $T^\ell(\Gamma',y',\beg)=(\Gamma,y,\nd)$ for some integer $\ell$ and frog $y\in Y(a,b)$.
	Thanks to the fact that $R$ is an involution and \Cref{claim:reverse}, we observe that
	\[
		RT^{\ell}R(\Gamma,y,\nd)=RT^{\ell}RT^{\ell}(\Gamma',y',\beg)=(\Gamma',y',\beg).
	\]
	Hence, we may recover $(\Gamma',y')$ from $(\Gamma,y)$. Thus, the map sending $(\Gamma',y')$ to $(\Gamma,y)$ is injective.
	This is an injection from the in-edges of $S$ to the frogs in $S$, so $\deg^{\text{in}}(S)\leq a+b$ for every $S\in {\Z/k\Z\choose a}\times{\Z/k\Z\choose b}$.
	We conclude that $\deg^{\text{in}}(S)=a+b$ for every $S\in {\Z/k\Z\choose a}\times{\Z/k\Z\choose b}$ via the hand-shaking lemma.
\end{proof}

Combining the above lemmas, we conclude:
\begin{theorem}
	For any $1\leq a,b\leq k$, the chain $P(a,b)$ admits a stationary distribution which is uniform on ${\Z/k\Z\choose a}\times{\Z/k\Z\choose b}$.
\end{theorem}
\begin{proof}
	Let $G_1,\dots,G_\ell$ be the weakly connected components of $G$.
	\Cref{lem:regular} shows that $G$ is a regular digraph; hence each $G_i$ is strongly connected.
	Since each $G_i$ is aperiodic (\Cref{lem:aperiod}), the random walk on any $G_i$ admits a unique stationary distribution $\pi_i$.
	Because $G_i$ is regular (\Cref{lem:regular}), the distribution is uniform on $G_i$, i.e.\ $\pi_i(S)=1/|V(G_i)|$ for each $S\in V(G_i)$ and $\pi_i(S)=0$ for each $S\notin V(G_i)$.

	The stationary distributions on $P(a,b)$ are precisely the convex combinations of $\pi_1,\dots,\pi_\ell$.
	In particular, the convex combination ${1\over |V(G)|}\bigl(|V(G_1)|\pi_1+\dots+|V(G_\ell)|\pi_\ell\bigr)$ is a stationary distribution of $P(a,b)$, which is uniform.
\end{proof}

Using our knowledge of the stationary distributions of $\hat P(m+1,m)$ and the coupling of this chain with the frog dynamics, we can finally prove \Cref{thm:margins}.

\begin{proof}[Proof of \Cref{thm:margins}]
	Let $\pi$ be the stationary distribution of the frog dynamics associated with $W=12\cdots k$ and let $F\sim\pi$.
	Due to \Cref{thm:nicespeed}, we know that $\Pr[F(\froggie_{\leq m})=T]={k\choose m}^{-1}$ for any $m\in[k]$ and $T\in{\{\lilypad_0,\dots,\lilypad_{k-1}\}\choose m}$.
	Therefore, by coupling $\hat P(m+1,m)$ and the frog dynamics, we compute
	\begin{align*}
		\Pr\bigl[F(\froggie_{m+1})=\lilypad_{\ell_{m+1}}\ \big|\ F(\froggie_{\leq m})=\{\lilypad_{\ell_1},\dots,\lilypad_{\ell_m}\}\bigr] &={k\choose m}\Pr\bigl[ f(S)=\bigl(\lilypad_{\ell_{m+1}},\{\lilypad_{\ell_1},\dots,\lilypad_{\ell_m}\}\bigr)\bigr]
	\end{align*}
	where $S$ is uniformly distributed on ${\Z/k\Z\choose m+1}\times{\Z/k\Z\choose m}$ since $P(m+1,m)$, and thus $\hat P(m+1,m)$, admits a uniform stationary distribution.
	Since $S$ is uniformly distributed, the conditional probability above is equal to
	\begin{equation}\label{eqn:dyck}
	{1\over{k\choose m+1}}\biggl|\biggl\{S\in{\Z/k\Z\choose m+1}\times{\Z/k\Z\choose m}:f(S)=\bigl(\lilypad_{\ell_{m+1}},\{\lilypad_{\ell_1},\dots,\lilypad_{\ell_m}\}\bigr)\biggr\}\biggr|.
	\end{equation}
	In other words, with $S^-=\{\ell_1,\dots,\ell_m\}$ fixed, we need to count the number of ways to select $S^+\in{\Z/k\Z\choose m+1}$ such that the optimistic frog occupies lily pad $\ell_{m+1}$.
	In order for the optimistic frog to occupy lily pad $\ell_{m+1}$, it must be the case that
	\[
		\sum_{i=\ell_{m+1}}^{\ell_{m+1}+j}\bigl(\mathbf 1[i\in S^+]-\mathbf 1[i\in S^-]\bigr)>0
	\]
	for all $j\in\{0,\dots,k-1\}$.

	Set $\beta_{m-i}=\bigl|S^+\cap(\ell_{m-i+1},\ell_{m-i}]\bigr|$ for $i\in\{0,\dots,m-1\}$.
	Rephrasing the above requirement, $S^+$ is valid if and only if $\sum_{i=0}^j\beta_{m-i}\geq j+1$ for all $j\in\{0,\dots,m-1\}$ and $\sum_{i=0}^{m-1}\beta_{m-i}=m$.
	Thus, setting $\Delta_i=\ell_i-\ell_{i+1}$ for $i\in[m]$, and making the substitution $\beta_{m-i}=\alpha_i$, \eqref{eqn:dyck} becomes
	\begin{align*}
		{1\over {k\choose m+1}}\sum_{\substack{\beta_m,\dots,\beta_1\geq 0\\ \sum_{i=0}^j\beta_{m-i}\geq j+1\\ \sum_{i=0}^{m-1}\beta_{m-i}=m}}\prod_{i=1}^m{\Delta_i\choose\beta_i}&={1\over{k\choose m+1}}\sum_{\substack{\alpha_1,\dots,\alpha_m\geq 0\\ \sum_{i\leq j}\alpha_i\leq j\\ \sum_{i\leq m}\alpha_i=m}}\prod_{i=1}^m{\Delta_i\choose\alpha_i}.\qedhere
	\end{align*}
\end{proof}

\section{Computer simulations}\label{sec:computerjoin}
\subsection{LCS between two random words}
\paragraph{New algorithm.}
Because of how easy it is to compute the LCS between a pair of words by using standard dynamic programming techniques, many researchers have computed the LCS between
random words (c.f.~\cite{lcsforests,chvatal_1975,liu_houdre_simulations}). They all used a fast $\widetilde{O}(n^2)$ time deterministic algorithm
to compute the LCS between words of length about $n$; this allowed them to perform extensive simulations. In contrast,
we used a faster $\widetilde{O}(n^{3/2})$ time probabilistic algorithm, which we discuss below. This new algorithm is not meant
to compute the LCS between two arbitrary words; it is designed to operate correctly only if its input is a pair of random
words. Unfortunately, unlike the slower algorithms used in the past, we have no rigorous guarantees that our algorithm actually
produces the correct output. Instead, we have only Monte Carlo simulations that compare this algorithm against a known correct
algorithm on the same set of inputs.

The basis for the algorithm is a conjecture \cite{houdre_diagonal} that an LCS between two random words
of equal length one should not match faraway symbols. More precisely, assume that $V,W\sim \bet^n$ are
two random words. Suppose $\LCS(V,W)=\ell$, and let $(V_{i_1}V_{i_2}\dotsb V_{i_\ell},\ W_{i_1}W_{i_2}\dotsb W_{i_\ell})$ be
a common subsequence of length~$\ell$. It is then conjectured that $\abs{i_r-j_r}$ can never be too large.
Geometrically, the conjecture says that the $\ell$ points $(i_1,j_1),\dotsc,(i_\ell,j_\ell)$ should be close
to the diagonal line $y=x$. The only rigorous result that we are aware of is that the points are asymptotically confined, with high probability,
to the wedge between the lines $y=c_1x$ and $y=c_2x$, for constants $c_1$, $c_2$ depending on the alphabet size \cite{houdre_diagonal}.

It is widely believe that the variance of $\LCS(V,W)$ is linear\footnote{It follows from Azuma's inequality that $\Var \LCS(V,W)$ is at most linear. However, it is not
known that $\Var \LCS(V,W)$ evens tends to infinity with $n$.};
it is thus reasonable to conjecture that most pairs $(i_r,j_r)$ satisfy $\abs{i_r-j_r}\leq C\sqrt{n}$ if $C$ is large. This suggests trying
to find a long common subsequence between $V$ and $W$ by restricting to only subsequences satisfying $\abs{i_r-j_r}\leq T$ for suitable $T$.
Using dynamic programming, this can done in time $O(nT)$.

In our implementation, we let $T_0=\lfloor\sqrt{2n}\rfloor$, and $T_i=\lfloor \tfrac{5}{2}T_{i-1}\rfloor$, and then used the method above for values of $T$ equal to $T_0,T_1,T_2,\dotsc$ in order until two
consecutive computations produced identical answers; that common value is the output value of our algorithm. The constant $\tfrac{5}{2}$ was chosen by accident; we did not try to optimize it.
The actual code used is available at \url{http://www.borisbukh.org/code/lcsfrogs19.html} and also as an ancillary file in the arXiv version of this paper.

To verify the algorithm's correctness, we ran a number of simulations in which we used both the new algorithm and the usual dynamic programming
algorithm on many pairs of random words $V,W\in \bet^n$. For the binary alphabet, we ran the experiment
with $n=2,500$ (300,000 times),
$n=5,000$   (300,000 times),
$n=10,000$  (100,000 times),
$n=20,000$  (30,000 times),
$n=50,000$  (18,350 times),
$n=250,000$ (2,670 times),
$n=500,000$ (1,080 times), and
$n=1,200,000$ (120 times).
In all these experiments, the new probabilistic algorithm produced the same answer as the usual dynamic programming algorithm.

These experiments gave us enough statistical evidence to trust the results of further computations using this algorithm.

After writing this paper, we learned from Alex Tiskin that a similar idea was independently proposed by Schimd and Bilard~\cite{schimd_edit} in the context of Levenshtein distance.
The \emph{Levenshtein distance} between words $V$ and $W$ is the least number of insertions, deletions, or symbol alterations needed
to turn $V$ into $W$. Note that the least number of insertions and deletions need to turn $V$ into $W$ is $\len V+\len W-2\LCS(V,W)$.
Denoting by $L(V,W)$ the Levenshtein distance between words $V$ and $W$, one may define $L_T(V,W)$ in a manner similar to how $\LCS_T(V,W)$ was defined.
Schimd and Bilardi computed $L_T$ for $T=\sqrt n$ and $n=2^{20}$ as a way to estimate $\E_{V,W\sim\bet^n}L(V,W)$.
To test their approach in the context of the LCS, we computed $\LCS_{\sqrt n}$ for $n=2^{20}$.
Out of 130 random trials, not in a single trial did $\LCS(V,W)$ and $\LCS_{\sqrt n}(V,W)$ agree.

\paragraph{Computational results and new conjectures.} Since the faster algorithm has allowed us to perform
more extensive (and thus more accurate) computations than before, these computations suggested new conjectures about
the LCS of a pair of random words.

To introduce the most interesting conjecture, consider a random word $V$ of length $2n$ as a concatenation of two
random words $V_1,V_2$, each of length $n$. Similarly, consider a random word $W$ of length $2n$ as a concatenation
of $W_1$ and $W_2$. It is then clear that $\LCS(V,W)\geq \LCS(V_1,W_1)+\LCS(V_2,W_2)$. From \cite{alexander_1994} one may deduce
that
\[
	\Delta(V,W)\eqdef \LCS(V,W)-\LCS(V_1,W_1)-\LCS(V_2,W_2)
\]
satisfies $\E\Delta(V,W)=O(\sqrt{n \log n})$.
Define $\Delta(2n)\eqdef\Delta(V,W)$ where $V,W$ are independent, random words of length $2n$.
Computing $\Delta(n)$ experimentally suggests the following.
\begin{conj}\label{cuberootconj}
	There are constants $c_1,c_2$ such that $\E\Delta(n)\sim c_1n^{1/3}$ and $\sqrt{\Var\Delta(n)}\sim c_2n^{1/3}$.
\end{conj}
The computational data behind the conjecture for the binary alphabet is summarized below.
\begin{center}
	\begin{tabular}{r|r|r|r}
		$n$ & number of trials& $\E\Delta(n)$ & $\sqrt{\Var\Delta(n)}$\\
		\hline\hline
		5,000 & 2102122 & 7.34957&4.41726\\
		10,000 & 3373157 & 9.46013&5.56865\\
		20,000 & 3225713 & 12.1248&7.01030\\
		40,000 & 505844 & 15.4730&8.81207\\
		80,000 & 68837 & 19.7529&11.1599\\
		160,000 & 40136 & 25.1560&14.0003\\
		320,000 & 95817 & 31.7925&17.6049\\
		640,000 & 19937 & 40.2075&22.0874\\
		1,280,000 & 10245 & 50.4519&27.8588\\
		2,560,000&7715&64.5401&34.6783\\
		5,120,000&1140&81.4482&44.5223
	\end{tabular}
\end{center}
In particular, the data suggests that $c_1\approx 1/2$ and $c_2\approx 1/4$ for the binary alphabet.
It is likely that $\Delta(n)/n^{1/3}$ converges to a non-trivial distribution; we do not have a conjecture
as to what that distribution is.

By summing $\E 2^{-i}\Delta(2^in)$ for $i=1,2,\dotsc$, we see that
the conjecture implies that
\[
  \E \LCS(V,W)=\gamma n-c_1'n^{1/3}+o(n^{1/3})\qquad\text{if }V,W\sim \bet^n,
\]
where $c_1'=c_1/(\sqrt[3]{4}-1)$.

The conjecture also strongly suggests that $f(n)\eqdef \Var \LCS(V,W)$ should grow linearly with~$n$.
Indeed, if $X,Y$ are any two mean-zero random variables, then $\Var[X+Y]=\Var[X]+2\E[XY]+\Var[Y]\geq
\Var[X] +\Var[Y]-2\sqrt{\Var[X] \Var[Y]}$. Because $\LCS(V_1,W_1)$ and $\LCS(V_2,W_2)$ are independent,
the variance of $\LCS(V_1,W_1)+\LCS(V_2,W_2)$ is $2f(n)$, so
$f(2n)\geq 2f(n)-2c_2n^{1/3}\sqrt{2f(n)}$. From this one may deduce that if we find a single $n_0$ for which
$f(n_0)> c_2'n_0^{2/3}$, where $c_2'=4c_2^2/(2-\sqrt[3]{4})^2$, then
$f(n)=\Omega(n)$ for $n$ of the form $2^i n_0$.

Assuming \Cref{cuberootconj}, one can thus make a more refined guess for the Chv\'atal--Sankoff constant $\gamma$ from \eqref{eq:chvatal_sankoff}
using simulations for several values of $n$. We obtained $\gamma\approx 0.8122$ for the binary alphabet. This is higher than
the previous guess of $0.8118$ from~\cite[Table~2]{lcsforests}, lower than the previous guess of $0.8126$ from~\cite[Table~1]{bundschuh2001} and is inside the interval $(0.8120,0.8125)$ suggested in \cite[Section~2.4]{dancik_phd}.

\paragraph{Periodic words.} We implemented the algorithm to compute the leading-term constant $\gamma_W=\gamma_W(1)$ in the formula
for $\E\LCS(R,W^{(n)})$ in \Cref{thm:basic}. The code is available at \url{http://www.borisbukh.org/code/lcsfrogs19.html}
and also as an ancillary file in the arXiv version.

Interestingly, there appear to exist periodic words that are more similar to the random word than the random word is!
More precisely, we found periodic binary words for which the leading constant $\gamma_W(1)$ exceeds $0.8122$, which is our conjectured value
of~$\gamma$. The binary word with the largest $\gamma_W(1)$ that we found is $W=0110111010010110010001011010$, for which $\gamma_W(1)\geq 0.82118$.
Alas, we cannot prove that $\gamma_W(1)>\gamma$ since the best rigorous upper bound on $\gamma$ for the binary alphabet is $\gamma\leq 0.826280$ \cite{lueker_2009}.

\section{Remarks}\label{sec:remarks}
\begin{itemize}[leftmargin=*]
	\item We are mystified by the coupling used to prove the stationary distribution of the frog dynamics associated with $W=12\cdots k$.
		We found it by first guessing the formula in \Cref{thm:margins}, noticing its combinatorial interpretation as a count of certain Dyck paths, and then looking for a suitable coupling.
		However, we do not have any high-level explanation for the appearance of Dyck paths nor for the time-reversal property in \Cref{claim:reverse}.
        As pointed by one of the referees, the coupling is similar to that by Angel \cite{angel} for the multi-species TASEP on a ring, which might help find the explanation (see also \cite{ferrari_martin} for a queuing interpretation of Angel's coupling).
	\item There is a fast algorithm to compute the LCS between a periodic word and any other word.
		Indeed, \Cref{prop:hit} implies that for any $F\in \cal F$ and $a\in\bet$, we can compute $Fa$ and $D_1(F,a),\dots,D_k(F,a)$ in $O(k)$ operations.
		Thus, for any $R\in \bet^n$, we can compute $D_1(F_\varnothing,R),\dots,D_k(F_\varnothing,R)$ in $O(kn)$ operations.
		\Cref{thm:ledgeloc} then tells us that
		\[
			\LCS(R,W^{(x)})=h_R(x)=x-\sum_{i:\ x_i\leq x}\Bigl\lceil{x-x_i\over k}\Bigr\rceil,
		\]
		where $x_i=D_i(F_\varnothing,R)+i-1$.
		Thus, we can compute $\LCS(R,W^{(x)})$ in $O(kn)$ operations.

		An $O(kn)$-time algorithm of a similar flavor was given by Tiskin~\cite{tiskin_2009}.
	\item We gave an algorithm to compute $\gamma_W$ from $W$ which relies on computing $s_1,\dots,s_k$ from the stationary distribution of the auxiliary frog dynamics.
		The set of all frog arrangements has size $k!$ and thus the auxiliary frog dynamics has $|\bet|\cdot k!$ states.

		However, using the ideas in \Cref{sec:linear}, we can actually compute $s_1,\dots,s_k$ from the stationary distributions of much smaller chains.
		Indeed, the $m$-arrangement chain associated with the word $W$ has only ${k\choose m}$ states and thus, through an extension of \Cref{lem:unstrat} to arbitrary words, $\sum_{i=1}^m s_i$ can be computed from the stationary distribution of a chain on $|\bet|\cdot{k\choose m}$ states.
		This observation allows us to compute $\gamma_W$ by instead finding the stationary distributions of a chain with only $|\bet|\cdot 2^k$ states.
    \item What can be said about the LCS between a random word and an arbitrary fixed word?
        In particular, we are interested in the following question:

        Define the constants
        \[
            c_k\eqdef\lim_{n\to\infty}{1\over n}\max_{W\in[k]^n}\E_{R\sim[k]^n}\LCS(R,W).
        \]
        How much larger is $c_k$ compared to the Chv\'atal--Sankoff constant over the same alphabet?
        In \Cref{sec:computerjoin}, we showed that $c_2\geq 0.8211$, whereas we believe the Chv\'atal--Sankoff constant for the binary alphabet to be approximately $0.8122$.
        By the work of Kiwi--Loebl--Matou\v{s}ek~\cite{kiwi_2005} the Chv\'atal--Sankoff constant over a $k$-letter alphabet is asymptotic to $2/\sqrt{k}$.
        In the opposite direction, a straightforward application of the union bound yields $c_k\lesssim e/\sqrt k$.
        We suspect that $c_k\sim 2/\sqrt k$.
\end{itemize}

\bibliographystyle{plain}
\bibliography{lcsfrogs}

\appendix

\section{Markov chain central limit theorem}\label{apx:centlim}

Here we give a derivation of the following statement which was used in the proof of \Cref{lem:norm}.

Consider a Markov chain $X_0,X_1,\dots$ on a finite state-space $\Omega$ with a unique stationary distribution $\pi$ and let $f\colon\Omega\to\R$ be any function.
Let $G$ be the digraph describing the Markov chain, i.e., the digraph on the vertex set $\Omega$ with an edge $u\to v$ whenever the transition probability
from $u$ to $v$ is positive. Let $\Omega^*$ denote the support of $\pi$ and set $G^*\eqdef G[\Omega^*]$.
Note that $G^*$ is strongly connected.

\begin{theorem}
	If $G^*$ contains closed walks $u_0\to u_1\to\dotsc\to u_m=u_0$ and $v_0\to v_1\to\dotsc\to v_n=v_0$ with ${1\over m}\sum_{t=1}^m f(u_t)\neq{1\over n}\sum_{t=1}^nf(v_t)$, then ${1\over\sqrt{n}}\sum_{t=1}^n\bigl(f(X_t)-\E_{X\sim\pi}f(X)\bigr)\tod\cal N(0,\sigma^2)$ for some fixed $\sigma\neq 0$.
\end{theorem}
\begin{proof}
	Starting with $\tau_0\eqdef\min\{n>0:X_n=u_0\}$, define $\tau_k\eqdef\min\{n>\tau_{k-1}:X_n=u_0\}$.
	Set $\zeta_k\eqdef\tau_k-\tau_{k-1}$ and $Y_k\eqdef\sum_{t=\tau_{k-1}+1}^{\tau_k}f(X_t)$, so $\zeta_1,\zeta_2,\dots$ are i.i.d., as are $Y_1,Y_2,\dots$.
	Since $\pi(u_0)>0$, it follows from~\cite[Section 16, Theorem 1]{Ccentlim} that, denoting $M\eqdef \E Y_k/\E\zeta_k$ and $\bar\sigma^2\eqdef\E\bigl[(Y_k-\zeta_kM)^2\bigr]$, we have
	\[
		{1\over\sqrt{n}}\sum_{t=1}^n\bigl(f(X_t)-\E_{X\sim\pi}f(X)\bigr)\tod\cal N(0,\pi(u_0)\bar\sigma^2),
	\]
	provided that $0<\bar\sigma^2<\infty$.
	Certainly $\bar\sigma^2<\infty$, so we need argue only that $\bar\sigma\neq 0$, which amounts to arguing that $Y_k-\zeta_k M$ is not identically $0$.

	If it were to be the case that $Y_k-\zeta_kM\equiv 0$, then any closed walk $w_0\to w_1\to\dotsc\to w_\ell=w_0$ with $w_0=u_0$ must satisfy
	\[
		{1\over\ell}\sum_{t=1}^\ell f(w_t)=M.
	\]
	In particular, ${1\over\ell}\sum_{t=1}^\ell f(w_t)={1\over m}\sum_{t=1}^m f(u_t)$.

	Now, since $v_0\in\Omega^*$ as well, we can find a pair of walks $u_0=r_0\to r_1\to\dots\to r_{n_1}=v_0$ and $v_0=r_0'\to r_1'\to\dots\to r_{n_2}'=u_0$ in $G^*$.
	Consider the closed walk which starts at $u_0$, traverses $r_0\to\dots\to r_{n_1}$, moves around $v_0\to\dots\to v_n$ a total of $K$ times, and finally traverses $r_0'\to\dots\to r_{n_2}'$ back to $u_0$.
	By the observation above, we must have
	\[
		{1\over m}\sum_{t=1}^mf(u_t)={1\over n_1+n_2+Kn}\biggl(\sum_{t=1}^{n_1}f(r_t)+\sum_{t=1}^{n_2}f(r_t')+K\sum_{t=1}^nf(v_t)\biggr),
	\]
	for every positive integer $K$.
	However, as $K\to\infty$, the right-hand side converges to ${1\over n}\sum_{t=1}^n f(v_t)$; contradicting our original assumption.
\end{proof}

\section{The maximum of two Gaussians}\label{apx:notgaus}

The following proposition, which was used in the proof of \cref{nnn:hitspeed} of \Cref{thm:normnotnorm}, was communicated to us by Tomasz Tkocz.
\begin{proposition}\label{prop:notgaus}
	If $X$ and $Y$ are centered (possibly degenerate) Gaussian random variables with $\Pr[X=Y]<1$, then $\max\{X,Y\}$ is \emph{not} a Gaussian random variable.
\end{proposition}
\begin{proof}
	Suppose that $X\sim\cal N(0,a^2)$ and $Y\sim\cal N(0,b^2)$ and set $Z\eqdef\max\{X,Y\}={1\over 2}\bigl(X+Y+|X-Y|\bigr)$.
	Suppose for the sake of contradiction that $Z\sim\cal N(\mu,\sigma^2)$ for some $\mu,\sigma$.
	Consider the moment generating function of $Z$: $\E e^{tZ}=e^{t^2\sigma^2/2+t\mu}$ for $t\in\R$.
	For $t>0$, we have the point-wise bounds,
	\[
		{e^{tX}+e^{tY}\over 2}\leq e^{tZ}\leq e^{tX}+e^{tY}\implies{e^{t^2a^2/2}+e^{t^2b^2/2}\over 2}\leq e^{t^2\sigma^2/2+t\mu}\leq e^{t^2a^2/2}+e^{t^2b^2/2}.
	\]

	If $\sigma^2>\max\{a^2,b^2\}$, then the inequality $e^{t^2\sigma^2/2+t\mu}\leq e^{t^2a^2/2}+e^{t^2b^2/2}$ is violated for sufficiently large $t$.

	If $\sigma^2<a^2$ or $\sigma^2<b^2$, then the inequality $e^{t^2\sigma^2/2+t\mu}\geq{1\over 2}\big(e^{t^2a^2/2}+e^{t^2b^2/2}\big)$ is violated for sufficiently large $t$.

	Finally, if $\sigma^2=a^2=b^2$, then $e^{t\mu}\leq 2$ for all $t>0$.
	However, since $\Pr[X=Y]<1$, we have
	\[
		\mu=\E Z={1\over 2}\E|X-Y|>0,
	\]
	and so this is impossible.
\end{proof}

\section{Case check for \texorpdfstring{\Cref{claim:reverse} (time-reversal)}{Claim~\ref{claim:reverse}}}\label{apx:reverse}
The table below verifies that $RTRT(\Gamma,y,t)=(\Gamma,y,t)$ for every $(\Gamma,y,t)\in\Omega(a,b)\setminus\Omega_\nd(a,b)$.
As we observed in \Cref{claim:reverse}, in order to verify this fact, we need keep track only of lily pads $x$ and $x+1$ where $x=\Gamma(y)$.

Each row of the table follows the trajectory of some $(\Gamma,y,t)$ under the successive maps $T$, $RT$, $TRT$, $RTRT$.
To reduce the number of cases, we use $?^-$ to denote a negative frog that might or might not be present.
In effect, the use of $?^-$ hides two cases: one with $?^-$ replaced by $r^-$ for some $r$, and one in which $?^-$ does not appear at all.
Similarly, the notation $?^+$ refers to a positive frog that might or might not be present.
Note that $R$ maps $?^-$ to $?^+$ and vice versa.

\def\s{1.1}
\tikzstyle{lily}=[draw,circle,minimum size=\s cm]
\def\frognul#1{\begin{tikzpicture}\node[lily] at (0,0) {};\node at (0,2*\s/3) {{\tiny $\mathclap{#1}$}};\end{tikzpicture}}
\def\frogone#1#2{\begin{tikzpicture}\node[lily] at (0,0) {};\node at (0,2*\s/3) {{\tiny $\mathclap{#1}$}};\node at (\s/15,0) {$\mathclap{#2}$};\end{tikzpicture}}
\def\frogtwo#1#2#3{\begin{tikzpicture}\node[lily] at (0,0) {};\node at (0,2*\s/3) {{\tiny $\mathclap{#1}$}};\node at ($ (90:\s/5) + (\s/15,0) $) {$\mathclap{#2}$};\node at ($ (-90:\s/5) +(\s/15,0) $) {$\mathclap{#3}$};\end{tikzpicture}}
\def\frogthree#1#2#3#4{\begin{tikzpicture}\node[lily] at (0,0) {};\node at (0,2*\s/3) {{\tiny $\mathclap{#1}$}};\node at ($ (90:\s/4) + (\s/15,0) $) {$\mathclap{#2}$};\node at ($ (210:\s/4) + (\s/15,0) $) {$\mathclap{#3}$};\node at ($ (-30:\s/4) + (\s/15,0) $) {$\mathclap{#4}$};\end{tikzpicture}}

\begin{center}
	%{\footnotesize
	\begin{longtable}{c|c|c|c|c}
		$\Gamma$ \verb|\\| $(y,t)$ & $T$ & $RT$ & $TRT$ & $RTRT$\\ \nopagebreak\hline\hline

		% i^+ begin
		%1
		\frogtwo{x}{i^+}{j^-}\frogtwo{x+1}{?^+}{?^-} & \frogtwo{x}{i^+}{j^-}\frogtwo{x+1}{?^+}{?^-} & \frogtwo{k-x-1}{?^-}{?^+}\frogtwo{k-x}{i^-}{j^+} & \frogtwo{k-x-1}{?^-}{?^+}\frogtwo{k-x}{i^-}{j^+} & \frogtwo{x}{i^+}{j^-}\frogtwo{x+1}{?^+}{?^-}\\ \nopagebreak
		$(i^+,\beg)$ & $(j^-,\nd)$ & $(j^+,\beg)$ & $(i^-,\nd)$ & $(i^+,\beg)$ \\ \hline

		%2
		\frogone{x}{i^+}\frogone{x+1}{?^-} & \frognul{x}\frogtwo{x+1}{i^+}{?^-} & \frogtwo{k-x-1}{i^-}{?^+}\frognul{k-x} & \frogone{k-x-1}{?^+}\frogone{k-x}{i^-} & \frogone{x}{i^+}\frogone{x+1}{?^-}\\ \nopagebreak
		$(i^+,\beg)$ & $(i^+,\nd)$ & $(i^-,\beg)$ & $(i^-,\nd)$ & $(i^+,\beg)$ \\ \hline

		%3
		\frogone{x}{i^+}\frogtwo{x+1}{j^+}{?^-} & \frognul{x}\frogthree{x+1}{i^+}{j^+}{?^-} & \frogthree{k-x-1}{i^-}{j^-}{?^+}\frognul{k-x} & \frogtwo{k-x-1}{j^-}{?^+}\frogone{k-x}{i^-} & \frogone{x}{i^+}\frogtwo{x}{j^+}{?^-}\\ \nopagebreak
		$(i^+,\beg)$ & $(j^+,\trans)$ & $(i^-,\trans)$ & $(i^-,\nd)$ & $(i^+,\beg)$ \\ \hline

		%4-2\frogone{x}{i^+}\frogone{x+1}{j^-} & \frognul{x}\frogtwo{x+1}{i^+}{j^-} & \frogtwo{k-x-1}{i^-}{j^+}\frognul{k-x} & \frogone{k-x-1}{j^+}\frogone{k-x}{i^-} & \frogone{x}{i^+}\frogone{x+1}{j^-} \\ \nopagebreak
		%$(i^+,\beg)$ & $(i^+,\nd)$ & $(i^-,\beg)$ & $(i^-,\nd)$ & $(i^+,\beg)$ \\ \hline

		%5-3\frogone{x}{i^+}\frogtwo{x+1}{j^+}{\ell^-} & \frognul{x}\frogthree{x+1}{i^+}{j^+}{\ell^-} & \frogthree{k-x-1}{i^-}{j^-}{\ell^+}\frognul{k-x} & \frogtwo{k-x-1}{j^-}{\ell^+}\frogone{k-x}{i^-} & \frogone{x}{i^+}\frogtwo{x+1}{j^+}{\ell^-} \\ \nopagebreak
		%$(i^+,\beg)$ & $(j^+,\trans)$ & $(i^-,\trans)$ & $(i^-,\nd)$ & $(i^+,\beg)$ \\ \hline

		% i^- begin
		%6
		\frogtwo{x}{i^-}{?^+}\frognul{x+1} & \frogone{x}{?^+}\frogone{x+1}{i^-} & \frogone{k-x-1}{i^+}\frogone{k-x}{?^-} & \frognul{k-x-1}\frogtwo{k-x}{i^+}{?^-} & \frogtwo{x}{i^-}{?^+}\frognul{x+1} \\ \nopagebreak
		$(i^-,\beg)$ & $(i^-,\nd)$ & $(i^+,\beg)$ & $(i^+,\nd)$ & $(i^-,\beg)$ \\ \hline

		%7
		\frogtwo{x}{i^-}{?^+}\frogone{x+1}{j^+} & \frogone{x}{?^+}\frogtwo{x+1}{i^-}{j^+} & \frogtwo{k-x-1}{i^+}{j^-}\frogone{k-x}{?^-} & \frogone{k-x-1}{j^-}\frogtwo{k-x}{i^+}{?^-} & \frogtwo{x}{i^-}{?^+}\frogone{x+1}{j^+} \\ \nopagebreak
		$(i^-,\beg)$ & $(j^+,\trans)$ & $(i^+,\trans)$ & $(i^+,\nd)$ & $(i^-,\beg)$ \\ \hline\pagebreak

		%8
		$\Gamma$ \verb|\\| $(y,t)$ & $T$ & $RT$ & $TRT$ & $RTRT$\\ \nopagebreak\hline\hline
		\frogtwo{x}{i^-}{?^+}\frogtwo{x+1}{j^-}{?^+} & \frogone{x}{?^+}\frogthree{x+1}{i^-}{j^-}{?^+} & \frogthree{k-x-1}{i^+}{j^+}{?^-}\frogone{k-x}{?^-} & \frogtwo{k-x-1}{j^+}{?^-}\frogtwo{k-x}{i^+}{?^-} & \frogtwo{x}{i^-}{?^+}\frogtwo{x+1}{j^-}{?^+} \\ \nopagebreak
		$(i^-,\beg)$ & $(j^-,\trans)$ & $(i^+,\trans)$ & $(i^+,\nd)$ & $(i^-,\beg)$ \\ \hline

		\frogtwo{x}{i^+}{j^-}\frogone{x+1}{?^-} & \frogone{x}{j^-}\frogtwo{x+1}{i^+}{?^-} & \frogtwo{k-x-1}{i^-}{?^+}\frogone{k-x}{j^+} & \frogone{k-x-1}{?^+}\frogtwo{k-x}{i^-}{j^+} & \frogtwo{x}{i^+}{j^-}\frogone{x+1}{?^-} \\ \nopagebreak
		$(i^+,\trans)$ & $(i^+,\nd)$ & $(i^-,\beg)$ & $(j^+,\trans)$ & $(i^+,\trans)$ \\  \hline

		%15
		\frogtwo{x}{i^+}{j^-}\frogtwo{x+1}{\ell^+}{?^-} & \frogone{x}{j^-}\frogthree{x+1}{i^+}{\ell^+}{?^-} & \frogthree{k-x-1}{i^-}{\ell^-}{?^+}\frogone{k-x}{j^+} & \frogtwo{k-x-1}{\ell^-}{?^+}\frogtwo{k-x}{i^-}{j^+} & \frogtwo{x}{i^+}{j^-}\frogtwo{x+1}{\ell^+}{?^-}\\ \nopagebreak
		$(i^+,\trans)$ & $(\ell^+,\trans)$ & $(i^-,\trans)$ & $(j^+,\trans)$ & $(i^+,\trans)$ \\ \hline

		%16-14\frogtwo{x}{i^+}{j^-}\frogone{x+1}{\ell^-} & \frogone{x}{j^-}\frogtwo{x+1}{i^+}{\ell^-} & \frogtwo{k-x-1}{i^-}{\ell^+}\frogone{k-x}{j^+} & \frogone{k-x-1}{\ell^+}\frogtwo{k-x}{i^-}{j^+} & \frogtwo{x}{i^+}{j^-}\frogone{x+1}{\ell^-} \\ \nopagebreak
		%$(i^+,\trans)$ & $(i^+,\nd)$ & $(i^-,\beg)$ & $(j^+,\trans)$ & $(i^+,\trans)$ \\ \hline

		%17-15\frogtwo{x}{i^+}{j^-}\frogtwo{x+1}{\ell^+}{p^-} & \frogone{x}{j^-}\frogthree{x+1}{i^+}{\ell^+}{p^-} & \frogthree{k-x-1}{i^-}{\ell^-}{p^+}\frogone{k-x}{j^+} & \frogtwo{k-x-1}{\ell^-}{p^+}\frogtwo{k-x}{i^-}{j^+} & \frogtwo{x}{i^+}{j^-}\frogtwo{x+1}{\ell^+}{p^-}\\ \nopagebreak
		%$(i^+,\trans)$ & $(\ell^+,\trans)$ & $(i^-,\trans)$ & $(j^+,\trans)$ & $(i^+,\trans)$ \\ \hline

		%18
		\frogthree{x}{i^+}{j^+}{?^-}\frogone{x+1}{?^-} & \frogtwo{x}{j^+}{?^-}\frogtwo{x}{i^+}{?^-} & \frogtwo{k-x-1}{i^-}{?^+}\frogtwo{k-x}{j^-}{?^+} & \frogone{k-x-1}{?^+}\frogthree{k-x}{i^-}{j^-}{?^+} & \frogthree{x}{i^+}{j^+}{?^-}\frogone{x+1}{?^-} \\ \nopagebreak
		$(i^+,\trans)$ & $(i^+,\nd)$ & $(i^-,\beg)$ & $(j^-,\trans)$ & $(i^+,\trans)$ \\ \hline

		%19
		\frogthree{x}{i^+}{j^+}{?^-}\frogtwo{x+1}{\ell^+}{?^-} & \frogtwo{x}{j^+}{?^-}\frogthree{x+1}{i^+}{\ell^+}{?^-} & \frogthree{k-x-1}{i^-}{\ell^-}{?^+}\frogtwo{k-x}{j^-}{?^+} & \frogtwo{k-x-1}{\ell^-}{?^+}\frogthree{k-x}{i^-}{j^-}{?^+} & \frogthree{x}{i^+}{j^+}{?^-}\frogtwo{x+1}{\ell^+}{?^-} \\ \nopagebreak
		$(i^+,\trans)$ & $(\ell^+,\trans)$ & $(i^-,\trans)$ & $(j^-,\trans)$ & $(i^+,\trans)$ \\ \hline

		\frogthree{x}{i^-}{j^-}{?^+}\frognul{x+1} & \frogtwo{x}{j^-}{?^+}\frogone{x+1}{i^-} & \frogone{k-x-1}{i^+}\frogtwo{k-x}{j^+}{?^-} & \frognul{k-x-1}\frogthree{k-x}{i^+}{j^+}{?^-} & \frogthree{x}{i^-}{j^-}{?^+}\frognul{x+1}\\ \nopagebreak
		$(i^-,\trans)$ & $(i^-,\nd)$ & $(i^+,\beg)$ & $(j^+,\trans)$ & $(i^-,\trans)$ \\ \hline

		%27
		\frogthree{x}{i^-}{j^-}{?^+}\frogone{x+1}{\ell^+} & \frogtwo{x}{j^-}{?^+}\frogtwo{x+1}{i^-}{\ell^+} & \frogtwo{k-x-1}{i^+}{\ell^-}\frogtwo{k-x}{j^+}{?^-} & \frogone{k-x-1}{\ell^-}\frogthree{k-x}{i^+}{j^+}{?^-} & \frogthree{x}{i^-}{j^-}{?^+}\frogone{x+1}{\ell^+} \\ \nopagebreak
		$(i^-,\trans)$ & $(\ell^+,\trans)$ & $(i^+,\trans)$ & $(j^+,\trans)$ & $(i^-,\trans)$ \\ \hline

		%28
		\frogthree{x}{i^-}{j^-}{?^+}\frogtwo{x+1}{\ell^-}{?^+} & \frogtwo{x}{j^-}{?^+}\frogthree{x+1}{i^-}{\ell^-}{?^+} & \frogthree{k-x-1}{i^+}{\ell^+}{?^-}\frogtwo{k-x}{j^+}{?^-} & \frogtwo{k-x-1}{\ell^+}{?^-}\frogthree{k-x}{i^+}{j^+}{?^-} & \frogthree{x}{i^-}{j^-}{?^+}\frogtwo{x+1}{\ell^-}{?^+} \\ \nopagebreak
		$(i^-,\trans)$ & $(\ell^-,\trans)$ & $(i^+,\trans)$ & $(j^+,\trans)$ & $(i^-,\trans)$ \\ \hline

\end{longtable}
%}
\end{center}

\end{document}